\documentclass[a4paper,preprint,10pt,3p]{elsarticle}
\usepackage[utf8]{inputenc}
\usepackage[T1]{fontenc}
\usepackage{ae}
\usepackage[english]{babel}
\usepackage{amsmath,amssymb,amstext,amsfonts}
\usepackage{units}
\usepackage[squaren]{SIunits}
\usepackage{bm}
\usepackage{graphicx}
\usepackage{pgf}
\usepackage{pgfplots}
\usepackage{tikz}
\usetikzlibrary{arrows}
\usetikzlibrary{shapes}
\usetikzlibrary{decorations.text}
\usepackage{xcolor}
\usepackage{color}
\usepackage{tikzscale}
\usepackage{pgfplots}

\usetikzlibrary{shapes.geometric}
\usetikzlibrary{positioning}
\usetikzlibrary{positioning,shapes.multipart, fit,backgrounds,calc}

\usepackage{nicefrac}
\usepackage{array,tabularx}
\usepackage{multirow}
\newcolumntype{L}[1]{>{\raggedright\arraybackslash}p{#1}}
\newcolumntype{C}[1]{>{\centering\arraybackslash}p{#1}}
\newcolumntype{R}[1]{>{\raggedleft\arraybackslash}p{#1}}
\usepackage{booktabs}
\usepackage{mathtools}

\usepackage{geometry}
\usepackage{pdfpages}

\usepackage[colorinlistoftodos,prependcaption,textsize=tiny]{todonotes}
\usepackage{caption}
\usepackage{subcaption}
\usepackage{framed}
\usepackage[pdfencoding=auto]{hyperref}
\usepackage{placeins}

\usepackage[amsthm]{ntheorem}
\theoremheaderfont{\normalfont\bfseries}

\theoremstyle{remark}

\theoremstyle{break}

\newcommand{\mbR}{\mathbb{R}}
\newcommand{\mcO}{\mathcal{R}}
\newcommand{\bigo}[1]{\mcO\left(#1\right)}

\newcommand{\intd}{\textrm{\;d}}
\newcommand{\ddt}{\frac{\intd}{\intd t}}

\DeclareMathOperator{\average}{average}
\DeclareMathOperator{\rank}{rank}

\DeclareMathOperator{\Ima}{Im}

\theoremstyle{plain}

\newtheorem{mytheo}{Theorem}[section]

\theoremstyle{plain}

\newcommand{\hs}[1]{\hspace{#1 em}}
\newcommand{\h}{\hspace{1 em}}

\newcommand{\ybc}{{y_\mathrm{bc}}}
\newcommand{\tybc}{\tilde y_\mathrm{bc}}
\newcommand{\Nbc}{{N_\mathrm{bc}}}

\newcommand{\Rbc}{{R_\mathrm{bc}}}
\newcommand{\RbcII}{{R_\mathrm{bc}^2}}

\newcommand{\Xbc}{X_\mathrm{bc}}
\newcommand{\Xhom}{X_\mathrm{hom}}

\newcommand{\Vhhom}{V_{h,\mathrm{hom}}}
\newcommand{\Vhinhom}{V_{\mathrm{inhom}}}

\newcommand{\phihom}{\Phi_{\mathrm{hom}}}
\newcommand{\phiinhom}{\Phi_{\mathrm{inhom}}}
\newcommand{\phibc}{\Phi_{\mathrm{bc}}}

\newcommand{\Vrhom}{V_{r,\mathrm{hom}}}

\newcommand{\tVinhom}{\tilde{V}_{\mathrm{inhom}}}

\newcommand{\ahom}{a_{\mathrm{hom}}}

\newcommand{\abc}{a_{\mathrm{bc}}}

\newcommand{\Rinhom}{{R_{\mathrm{inhom}}}}
\newcommand{\Finhom}{F_{\mathrm{inhom}}}

\newcommand{\Rhom}{{R_{\mathrm{hom}}}}
\newcommand{\RhomII}{{R_{\mathrm{hom}}^2}}
\newcommand{\RhomIII}{{R_{\mathrm{hom}}^3}}

\newcommand{\Vinhom}{V_{\mathrm{inhom}}}

\newcommand{\Oh}{{\Omega_h}}
\newcommand{\Ohi}{{\Omega_h^{-1}}}

\newcommand{\tstart}{t_\mathrm{start}}
\newcommand{\tend}{t_\mathrm{end}}

\newcommand{\consistent}{\textit{varying angle}}
\newcommand{\Consistent}{Testcase \textit{varying angle}.}
\newcommand{\movemode}{\textit{moving mode}}
\newcommand{\Movemode}{Testcase \textit{moving mode}.}

\newcommand{\Vrvo}{V_{r,\mathrm{velocity-only}}}
\newcommand{\Vrvp}{V_{r,\mathrm{velocity-pressure}}}

\newcommand{\Li}{\bar L_h^{-1}}
\newcommand{\Lit}{\bar L_h^{-T}}
\newcommand{\bL}{\bar L_h}

\pgfplotsset{compat=newest}
\usepgfplotslibrary{groupplots}
\usepgfplotslibrary{polar}
\usepgfplotslibrary{smithchart}
\usepgfplotslibrary{statistics}
\usepgfplotslibrary{dateplot}
\usepgfplotslibrary{ternary}
\usetikzlibrary{arrows.meta}
\usetikzlibrary{backgrounds}
\usepgfplotslibrary{patchplots}
\usepgfplotslibrary{fillbetween}
\pgfplotsset{%
    layers/standard/.define layer set={%
            background,axis background,axis grid,axis ticks,axis lines,axis tick labels,pre main,main,axis descriptions,axis foreground%
        }{grid style= {/pgfplots/on layer=axis grid},%
            tick style= {/pgfplots/on layer=axis ticks},%
            axis line style= {/pgfplots/on layer=axis lines},%
            label style= {/pgfplots/on layer=axis descriptions},%
            legend style= {/pgfplots/on layer=axis descriptions},%
            title style= {/pgfplots/on layer=axis descriptions},%
            colorbar style= {/pgfplots/on layer=axis descriptions},%
            ticklabel style= {/pgfplots/on layer=axis tick labels},%
            axis background@ style={/pgfplots/on layer=axis background},%
            3d box foreground style={/pgfplots/on layer=axis foreground},%
        },
}

\graphicspath{{./figures/}{./figures/inkscape/}}

\setuptodonotes{fancyline}

\begin{document}

\begin{frontmatter}

% Title Page
% \title{A mass-conserving and energy-consistent reduced-order model for the incompressible Navier-Stokes equations with time-dependent boundary conditions}
\title{No pressure? %Mass-conserving and
Energy-consistent ROMs for the incompressible Navier-Stokes equations with time-dependent boundary conditions}
\author[1]{H. Rosenberger}
\author[1]{B. Sanderse}

% \author{H. Rosenberger}
% \author{B. Sanderse}

\address[1]{Centrum Wiskunde \& Informatica, Science Park 123, Amsterdam, The Netherlands}

% \affiliation[1]{organization={Centrum Wiskunde \& Informatica},%Department and Organization
%             addressline={Science Park 123}, 
%             city={Amsterdam},
%             country={The Netherlands}}

% \affiliation[2]{organization={Eindhoven University of Technology},%Department and Organization
%             addressline={De Zaale}, 
%             city={Eindhoven},
%             % postcode={5600 MB}, 
%             % state={State Two},
%             country={The Netherlands}}

\begin{abstract}
%% Text of abstract

This work presents a novel reduced-order model (ROM) for the incompressible Navier-Stokes equations with time-dependent boundary conditions.
This ROM 
is velocity-only, i.e.\ the simulation of the velocity does not require the computation of the pressure, and preserves the structure of the kinetic energy evolution.

The key ingredient of the novel ROM is a decomposition of the velocity into a field with homogeneous boundary conditions and a lifting function that satisfies the mass equation with the prescribed inhomogeneous boundary conditions. This decomposition is inspired by the Helmholtz-Hodge decomposition and exhibits orthogonality of the two components. 
% We have discovered this orthogonality to be crucial to preserve the structure of the kinetic energy evolution.
This orthogonality is crucial to preserve the structure of the kinetic energy evolution.
% For specific boundary conditions, this kinetic energy evolution implies a bound on the kinetic energy and hence nonlinear stability of the proposed ROM.
% Since classical hyper-reduction techniques are not applicable, we propose a novel method to make the ROM simulation efficient. This novel method involves the explicit approximation of the boundary condition contributions and hence benefits the interpretability of the ROM.
To make the evaluation of the lifting function efficient, we propose a novel method that involves an explicit approximation of the boundary conditions with POD modes, while preserving the orthogonality of the velocity decomposition and thus the structure of the kinetic energy evolution.

We show that the proposed velocity-only ROM is equivalent to a velocity-pressure ROM, i.e., a ROM that simulates both velocity and pressure. 
This equivalence can be generalized to other existing velocity-pressure ROMs and reveals valuable insights in their behaviour.

Numerical experiments on test cases with inflow-outflow boundary conditions confirm the correctness and efficiency of the new ROM, and the equivalence with the velocity-pressure formulation.

% thanks to the orthogonality of this decomposition, the ROM preserves the structure of the FOM kinetic energy evolution

% A new reduced order model (ROM) formulation for the incompressible Navier-Stokes equations is 
% proposed, which incorporates time-dependent boundary conditions while ensuring a divergence-free velocity field. As a result, the ROM is velocity-only, and can be simulated completely without knowledge of the pressure. The approach is to split the velocity into a homogeneous and an inhomogeneous component: the former is approximated with a standard POD-Galerkin setting, while for the latter we propose a new decomposition that is based on a POD of snapshots of the boundary values. 
% % For the considered numbers of modes, the speed-up of the ROM online phase compared to the FOM is approximately two orders of magnitude. 
% % In a numerical experiment, the ROM demonstrates a significantly improved accuracy compared to a standard POD-Galerkin ROM. 

% intro: 
% - existing ROM methods do not discuss arbitrary time-dependent boundary conditions
% - existing ROM lifting functions do not exhibit orthogonality which we have discovered to be a key property for preserving the structure of the kinetic energy evolution
\end{abstract}

\end{frontmatter}

\textbf{Keywords:} 
Incompressible Navier-Stokes equations; Reduced-order model; Time-dependent boundary conditions; Structure preservation; Energy consistency; POD-Galerkin

\section{Introduction}

% intro: 
% - existing ROM methods do not discuss arbitrary time-dependent boundary conditions
% - existing ROM lifting functions do not exhibit orthogonality which we have discovered to be a key property for preserving the structure of the kinetic energy evolution

Turbulent flow is an ubiquitous phenomenon in fluid flows such as 
blood flow in arteries, ocean currents and air flow around wind turbines. 
% The
Accurate
numerical simulation of these turbulent flows 
% is typically performed with 
requires
computationally 
expensive methods such as Direct Numerical Simulation (DNS) 
\cite{moin1998direct} or Large Eddy Simulation (LES) \cite{piomelli1999large}. 
The computational cost of these methods can sometimes be tolerated for a small number of simulations but becomes prohibitive when it comes to tasks that require many simulation runs such as in the case of optimization and uncertainty quantification studies. 

To mitigate the computational costs, several techniques have been developed to 
reduce the complexity of the respective full order model (FOM) resulting in 
% a
so-called reduced-order models (ROMs).
Many of these techniques are projection-based, i.e., they project the 
high-dimensional discretization of the model obtained from classical 
discretization methods such as Finite Elements, Finite Volumes or Finite 
Differences onto a low-dimensional space, the reduced basis space. Methods 
to obtain such reduced bases include sampling strategies 
\cite{haasdonk2008reduced} and 
Proper Orthogonal Decomposition (POD) \cite{volkwein2011model}. 

The POD, which is also known as Karhunen-Loève expansion and Principal 
Component Analysis, is optimal in the sense that it extracts the reduced basis 
from given 
high-resolution data which minimizes the approximation error with respect to 
this data in the $l^2$ norm \cite{volkwein2011model}. However, several authors 
observed that ROMs resulting from the Galerkin projection of such POD reduced 
bases tend 
to exhibit instabilities when applied to fluid flow 
problems \cite{fick2018stabilized,stabile2018finite}. This observation is often related 
to energy dissipation at the smallest scales of turbulent flow, which are not included 
in the POD reduced basis \cite{fick2018stabilized}. Hence, approaches to 
address the ROM 
instability include methods that add such small scale effects either by 
introducing additional dissipation via closure modeling \cite{wang2012proper,ahmed2021closures} or using a POD based 
on a discrete equivalent of the $H^1$ norm instead of the $l^2$\,norm in order 
to achieve a better 
representation of small scales in the reduced basis \cite{iollo2000stability}.

Another cause of instability has been identified in the ROM pressure field. 
This pressure field can violate an $\inf$-$\sup$ condition at the ROM level,
even if the $\inf$-$\sup$ condition of the underlying FOM is satisfied \cite{ballarin2015supremizer,stabile2018finite}.
When violating the condition, the ROM can exhibit spurious pressure modes that induce instability.
Hence, many works have been developed to satisfy this condition on the ROM level \cite{ballarin2015supremizer,stabile2018finite,gerner2012certified,caiazzo2014numerical,rozza2013reduced}.
% which is not ensured to be stable if an $\inf$-$\sup$ condition is violated 
% \cite{ballarin2015supremizer}. Hence, methods have been developed which satisfy 
% this condition, e.g., \cite{ballarin2015supremizer,stabile2018finite}.
% At the same time,
On the other hand,
it has been shown that even for a linear ordinary 
differential equation (ODE), a stable FOM 
can result in an unstable POD-Galerkin ROM \cite{prajna2003pod}. 
% As a result, 
Consequently, POD-Galerkin ROM instabilities are not exclusively 
related to the incompressible Navier Stokes equations. 

A promising approach to avoid ROM instabilities which is not limited to specific equations is the idea of 
structure preservation 
\cite{carlberg2018conservative,chan2020entropy,peng2016symplectic,sanderse2020non}.
 Based on the observation 
that instability is often observed in combination with a blow-up of 
certain physical quantities such as energy,
 Peng et al.\ \cite{peng2016symplectic} propose a
POD-Galerkin ROM that preserves the symplecticity of Hamiltonian systems. Similarly, Chan 
\cite{chan2020entropy} constructs ROMs for nonlinear conservation laws 
satisfying a semi-discrete entropy dissipation.
% Carlberg et al.\ \cite{carlberg2018conservative} improve the accuracy of 
% POD-Galerkin based ROMs for the Euler equations by enforcing global 
% conservation of all primary variables by solving constrained optimization 
% problems, and Schein et al.\ \cite{schein2021preserving} propose a constrained 
% optimization framework for general equality and inequality constraints for 
% structure preservation.
For the incompressible Navier-Stokes equations, the POD-Galerkin ROM proposed in \cite{sanderse2020non} achieves nonlinear stability by preserving the conservation of kinetic energy in the inviscid limit for homogeneous boundary conditions. In this work, 
% our aim is to 
% generalize this ROM 
this ROM is generalized
to arbitrary time-dependent inhomogeneous boundary conditions, while possessing three properties: (i) structure preservation (consistency with the kinetic energy equation), (ii) velocity-only formulation (pressure-free), and (iii) arbitrary inhomogeneous boundary conditions (not separable in space- and time-dependent functions). To our best knowledge such a ROM formulation does not yet exist.

To achieve the first two properties, we propose the use of a specific lifting function. Classically, lifting functions are chosen as a velocity field that satisfies the prescribed inhomogeneous boundary conditions, so that the ROM can be solved with homogeneous boundary conditions  \cite{mohebujjaman2017energy,grassle2019pod,fick2017reduced}.
Inspired by the Helmholtz-Hodge decomposition \cite{bhatia2012helmholtz}, we propose a lifting function that is orthogonal to the homogeneous velocity field. This orthogonality allows us to preserve the structure of the kinetic energy equation on the ROM level. 

% A popular approach to take inhomogeneous boundary conditions into account is the use of a lifting function \cite{mohebujjaman2017energy,grassle2019pod,fick2017reduced}.
% However, most of these approaches only consider time-independent 
% \cite{mohebujjaman2017energy} 
% or very simple time-dependent boundary conditions \cite{grassle2019pod}. In this work, in contrast, we propose a mass-conserving ROM for arbitrary time-dependent boundary conditions. 
% The lifting function of this ROM is inspired by the Helmholtz-Hodge decomposition \cite{bhatia2012helmholtz} on the FOM level.
% Such lifting functions are typically chosen as a velocity field that satisfies the prescribed inhomogeneous boundary conditions. By subtracting the lifting function from a solution to the incompressible Navier-Stokes equations with inhomogeneous boundary conditions, one obtains a velocity field that satisfies homogeneous boundary conditions. Existing approaches on how to construct such lifting functions include purely discretely-motivated node-wise functions \cite{grassle2019pod} and the solution to a Stokes problem \cite{mohebujjaman2017energy,fick2017reduced}.

% In contrast to the case of homogeneous boundary conditions in \cite{sanderse2020non}, preserving this structure does not directly imply nonlinear stability because the kinetic energy is in general not a conserved quantity for inhomogeneous boundary conditions. However, for specific boundary conditions, the kinetic energy is bounded and hence implies nonlinear ROM stability.

In addition, this choice of lifting function has the benefit that the pressure is eliminated from the formulation, leading to a velocity-only ROM, i.e., the ROM velocity can be simulated without requiring the computation of the pressure. An important benefit of eliminating the pressure is that our velocity-only ROM is not affected by the $\inf$-$\sup$ instability mentioned above.
We stress that the omission of the pressure is \textit{not} an approximation as performed in some early works \cite{ito1998reduced,peterson1989reduced}. As shown in \cite{noack2005need}, merely neglecting the pressure in ROMs can lead to significant errors. In our ROM, in contrast, the elimination of the pressure is rigorous and exact. In fact, we show that the proposed velocity-only ROM is equivalent to a velocity-pressure ROM (for a specific choice of basis). 
% Vice versa, we show that existing velocity-pressure ROMs can be interpreted as velocity-only ROMs. These equivalences give 
This equivalence gives us
valuable insights in the behaviour of the ROMs.

% addressed in many other works \cite{ballarin2015supremizer,stabile2018finite,gerner2012certified,caiazzo2014numerical,rozza2013reduced}.

Lastly, to be able to deal with arbitrary inhomogeneous time-dependent boundary conditions, we propose to approximate the boundary conditions by projection onto a POD basis. This significantly reduces the computational costs that are associated with dealing with the exact boundary conditions, and is more accurate and interpretable than the common way of using hyper-reduction \cite{chaturantabut2009discrete, brunton2019data}.

% While most other ROM methods only address homogeneous, time-independent inhomogeneous or time-dependent boundary conditions with a specific structure, our proposed ROM is 
% % in principle 
% applicable to arbitrary boundary conditions. However, dealing with the exact boundary conditions can be prohibitively time- or memory-consuming. A common way of mitigating such cost is hyper-reduction. In our ROM, however, hyper-reduction is not directly applicable. Instead, we approximate the boundary conditions by projection onto a POD basis. Thanks to this explicit boundary condition approximation our ROM has better interpretability than hyper-reduced ROMs.

% This article is divided into XX Sections.
This article is organized as follows.
In Section \ref{chap: intro to sanderse}, we introduce the continuous formulation and our finite volume discretization of the incompressible Navier-Stokes equations. In Section \ref{sec:velo only rom}, we present our novel velocity-only ROM for time-dependent boundary condition. In Section \ref{sec:kin energy evol}, we show that this ROM preserves the structure of the kinetic energy evolution. In Section \ref{sec:equivalent velo pres rom}, we show the ROM's equivalence to a velocity-pressure ROM. In Section \ref{sec:num ex}, we test the novel ROM in two test cases involving time-dependent inflow boundary conditions.

\section{Preliminaries: continuous equations, full-order model, and
% reduced-order model formulation}
projection-based model order reduction}
\label{chap: intro to sanderse}

\subsection{The incompressible Navier-Stokes equations 
% and Helmholtz-Hodge decomposition}
}

The incompressible Navier-Stokes equations describe the conservation of mass 
and momentum of a fluid in a domain $\Omega\subset\mbR^d$ over a time interval 
$[0,T]$, and can be expressed as 
\begin{align}
	\nabla \cdot u &= 0 &\textrm{in 
	}\Omega\times[0,T],
	\label{eq:conti mass cons} \\
	\frac{\partial}{\partial t} u + \nabla \cdot (u\otimes u) &= -\frac{1}{\rho} \nabla p + 
	\nu \nabla^2 u + f  &\textrm{in 
	}\Omega\times[0,T],
	\label{eq:conti mom cons}
	\intertext{with divergence-free initial condition}
	u(x,0) &= u_0(x) & \textrm{ in }&\Omega\times[0,T],
	\label{eq:hhd initial condition} 
\end{align}
where the vector field $u:\Omega\times[0,T]\rightarrow \mbR^d$ describes the velocity, the scalar field $p:\Omega \times [0,T] \rightarrow \mbR$ the pressure, $\nu$ is the kinematic viscosity, $\rho$ the (constant) density (which will be absorbed in the pressure in subsequent equations), and $f$ a forcing term (possibly space- and time-dependent). A set of boundary conditions which is appropriate to ensure uniqueness of the solution $u$ 
% is necessary 
(see \cite[Section 3.8]{gresho1998incompressible})
is assumed to be prescribed. 

\subsection{Full order model: Finite volume discretization}
\label{sec: a struc-pres fv disc}
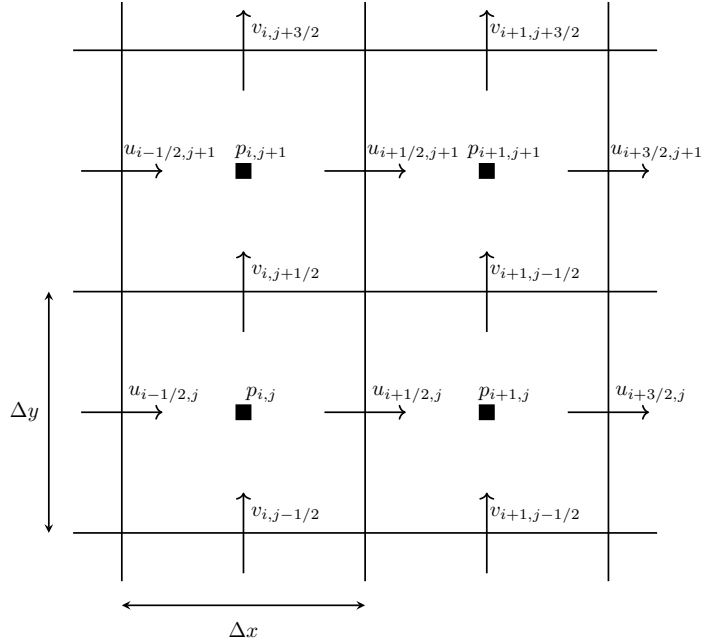
\begin{figure}[h]%[H] 
	\centering
	\scalebox{0.8}{	%\documentclass[tikz, margin=5mm]{standalone}
%\usepackage{tikz}
%\usetikzlibrary{shapes.geometric}
%\usetikzlibrary{positioning}
%\usetikzlibrary{positioning,shapes.multipart, fit,backgrounds,calc}
%
%%\tikzset{dot/.style={fill=black,circle}}
%%\tikzset{edot/.style={fill=black,circle}}
%
%\begin{document}
\begin{tikzpicture}[thick,black]
\tikzset{edot/.style={fill=black,diamond}}

\tikzset{triangle/.style = {fill=black!0, regular polygon, regular polygon 
sides=3}}
\tikzset{grey_triangle/.style = {fill=black!30, regular polygon, regular 
polygon sides=3}}
\tikzset{border rotated/.style = {shape border rotate=270}}
\tikzset{square/.style = {fill=black!100, rectangle}}

\foreach \x in {4}
{
	\foreach \y in {0,1,2}
	{
		\draw (-0.2*\x,\y*\x) -- (2.2*\x,\y*\x);
		\draw (\y*\x,-0.2*\x) -- (\y*\x,2.2*\x);
		
		\foreach \z in {0,...,1}
		{
			%\node [name=A] at (0.75*\x,0.5*\x) {};
			%\node [name=B] at (1.25*\x,0.5*\x) {};
			%\draw [->] (A) edge (B);
			\node [name=A] at (\y*\x-0.2*\x,\z*\x+0.5*\x) {};
			\node [name=B] at (\y*\x+0.2*\x,\z*\x+0.5*\x) {};	
			\draw [->] (A) edge (B);
% 			\node at  (\y*\x  +.3,\z*\x+0.5*\x  +.3) {$u$};
			
			\node [name=A] at (\z*\x+0.5*\x,\y*\x-0.2*\x) {};
			\node [name=B] at (\z*\x+0.5*\x,\y*\x+0.2*\x) {};	
			\draw [->] (A) edge (B);
% 			\node at (\z*\x+0.5*\x   +.3,\y*\x   +.3){$v$};	
		}
	}
	
% 	\draw [stealth-stealth] (2*\x,1.8*\x) -- (3*\x,1.8*\x);
% 	\draw [stealth-stealth] (2.4*\x,0*\x) -- (2.4*\x,1*\x);
% 	\node at (2.5*\x,1.9*\x) {$\Delta x$};
% 	\node at (2.6*\x,.5*\x) {$\Delta y$};

	\draw [stealth-stealth] (0*\x,-.3*\x) -- (1*\x,-.3*\x);
	\draw [stealth-stealth] (-.3*\x,0*\x) -- (-.3*\x,1*\x);
	\node at (.5*\x,-.4*\x) {$\Delta x$};
	\node at (-.4*\x,.5*\x) {$\Delta y$};
	
	\node at  (0*\x  +.7,0*\x+0.5*\x  +.3) {$u_{i-1/2,j}$};
	\node at  (0*\x  +.8,1*\x+0.5*\x  +.3) {$u_{i-1/2,j+1}$};
	\node at  (1*\x  +.7,0*\x+0.5*\x  +.3) {$u_{i+1/2,j}$};
	\node at  (1*\x  +.8,1*\x+0.5*\x  +.3) {$u_{i+1/2,j+1}$};
	\node at  (2*\x  +.7,0*\x+0.5*\x  +.3) {$u_{i+3/2,j}$};
	\node at  (2*\x  +.8,1*\x+0.5*\x  +.3) {$u_{i+3/2,j+1}$};
	
	\node at (0*\x+0.5*\x   +.7,0*\x   +.3){$v_{i,j-1/2}$};
	\node at (1*\x+0.5*\x   +.8,0*\x   +.3){$v_{i+1,j-1/2}$};
	\node at (0*\x+0.5*\x   +.7,1*\x   +.3){$v_{i,j+1/2}$};
	\node at (1*\x+0.5*\x   +.8,1*\x   +.3){$v_{i+1,j-1/2}$};
	\node at (0*\x+0.5*\x   +.7,2*\x   +.3){$v_{i,j+3/2}$};
	\node at (1*\x+0.5*\x   +.8,2*\x   +.3){$v_{i+1,j+3/2}$};
	
	\node at (0*\x+0.5*\x +.3,0*\x+0.5*\x +.3) {$p_{i,j}$};
	\node at (0*\x+0.5*\x +.3,1*\x+0.5*\x +.3) {$p_{i,j+1}$};
	\node at (1*\x+0.5*\x +.3,0*\x+0.5*\x +.3) {$p_{i+1,j}$};
	\node at (1*\x+0.5*\x +.3,1*\x+0.5*\x +.3) {$p_{i+1,j+1}$};

	\foreach \y in {0,1}
	{
		\foreach \z in {0,1}
		{
%			\node[square,draw] at ((\y+.5)*\x,(\z+.5)*\x){};
			\node[square,draw] at (\y*\x+0.5*\x,\z*\x+0.5*\x) {};
% 			\node at (\y*\x+0.5*\x +.3 ,\z*\x+0.5*\x +.3 ) {$p$};
%			\node[square,draw] at (\y,\z) {};
		}
	}
%\draw (0,0) -- (2*\x,0);
%\draw (0,\x) -- (2*\x,\x);
%\draw (0,0) -- (0,\x);
%\draw (2*\x,0) -- (2*\x,\x);
%\draw [line width=2mm] (\x,0) -- (\x,\x);
%
%%%\node[rectangle,draw] at (4,4){};
%\node[square,draw] at (0.5*\x,0.5*\x){};
%\node[square,draw] at (1.5*\x,0.5*\x){};
%\node at (0.6*\x,0.6*\x) {$p_{I,j}$};
%\node at (1.6*\x,0.6*\x) {$p_{I+1,j}$};
%%
%\node [name=A] at (0.75*\x,0.5*\x) {};
%\node [name=B] at (1.25*\x,0.5*\x) {};
%\draw [->] (A) edge (B);
%\node at (1.25*\x,0.4*\x) {$u_{I+1/2,j}$};
}

\end{tikzpicture}
%\end{document}
 }
	% \scalebox{.8}{	\input{tikz/staggered grid} }
	\caption{Staggered grid. Velocity values are located on the faces of the 
	pressure finite volumes.}
	\label{fig: staggered grid}
\end{figure}

The Navier-Stokes equations \eqref{eq:conti mass 
cons}-\eqref{eq:conti mom cons} in $d=2$ dimensions are spatially discretized with 
the second-order finite volume method on a 
staggered grid proposed in 
% \cite{harlow1965numerical}
\cite{sanderse2020non}
and depicted in Fig.\ \ref{fig: staggered grid}. 
The extension to $d=3$ is straightforward, and all the methods described in the subsequent sections also apply to that case.

The $u$- and $v$-component of the velocity, and the pressure $p$ are discretized by 
$u_h(t)\in\mbR^{N_u},\,v_h(t)\in\mbR^{N_v}$ and $p_h(t)\in\mbR^{N_p}$, 
respectively, where $N_p = N_x\times N_y$ is the number of ``pressure finite 
volumes'' and $N_u, N_v \approx N_p$ 
% (the exact relation depending on the type of boundary conditions) 
are 
the numbers of ``velocity component finite volumes'' which depend on the type of 
boundary conditions. The velocity component vectors are aggregated in one 
velocity vector $V_h(t)\in\mbR^{N_V},\, N_V = N_u+N_v$. 
% The explicit time dependence of $V_h$ and $p_h$ and similar quantities will be left out if no confusion can arise. 
The spatial discretization of %the integral form 
of \eqref{eq:conti mass 
	cons}-\eqref{eq:conti mom cons} can be summarized as 
% 	\cite{sanderse2020non}
\begin{align}
	M_h V_h(t) &= y_M(t) 
	:= F_M \ybc(t)
	,
	\label{eq:spacedisc navier stokes mass}\\
	\Omega_h \ddt V_h(t) &=
% 	F_h\left(V_h(t),p_h(t),\ybc(t),f(t)\right).
    F_h^{CD}(V_h(t),\ybc(t)) - G_h p_h(t).
	\label{eq:spacedisc navier stokes momentum}
\end{align}
Here, $\Omega_h\in\mbR^{N_V\times N_V}$ is a diagonal matrix with the sizes of 
the velocity component finite volumes corresponding to the entries of $V_h$ on 
its diagonal. The matrix $M_h\in\mbR^{N_p\times N_V}$ describes the mass 
equation \eqref{eq:conti mass cons} such that each 
% entry of the vector $M_h V_h$ 
row of \eqref{eq:spacedisc navier stokes mass}
represents conservation of mass over the finite volume around a pressure point 
$p_{i,j}$ as depicted in Figure \ref{fig: staggered grid},
\begin{align}
	u_{i+1/2,j} \Delta y - u_{i-1/2,j} \Delta y +  v_{i,j+1/2} \Delta x - v_{i,j-1/2} \Delta x = 0.
\end{align}
% \begin{align}
% 	\bar u_{i+1/2,j}-\bar u_{i-1/2,j} + \bar v_{i,j+1/2} - \bar v_{i,j-1/2} = 0.
% \end{align}
% The notation $\bar \cdot$ means integration over a finite volume face, 
% e.g.,\ $\bar u_{i+1/2,j} = u_{i+1/2,j}\Delta y$ 
using the midpoint rule.
% The vectors $\ybc(t)\in\mbR^\Nbc$ and $y_M(t)\in\mbR^{N_p}$ represent the boundary conditions as described in more detail in \ref{app:discretiz of bc} and the vector $f(t)\in\mbR^{N_V}$ represents external forces.

The vector $\ybc(t) \in\mbR^\Nbc$ constitutes the time-dependent parametrization of inhomogeneous boundary conditions. In this work, we only consider time-dependence of inflow boundary conditions. Hence, $\ybc(t)$ comprises the prescribed Dirichlet velocity values at finite volume faces on inflow boundaries and $\Nbc$ is the number of such faces. However, the construction of $\ybc(t)$ can be generalized to test cases with other types of time-dependent boundary conditions.
Note that $\Nbc$ scales only by $\bigo{h^{-d+1}}$, while $N_V$ and $N_p$ scale by $\bigo{h^{-d}}$ upon decreasing the mesh size $h$.

The matrix $F_M\in\mbR^{N_p\times \Nbc}$ maps the vector $\ybc(t)$ to the mass equation right-hand side $y_M(t)\in\mbR^{N_p}$.
% , 
% and the vector $f(t)\in\mbR^{N_V}$ represents external forces.
% The matrix $F_M\in\mbR^{N_p\times \Nbc}$ maps the vector $\ybc(t)$, which represents the prescribed boundary conditions, to the mass equation right-hand side $y_M(t)\in\mbR^{N_p}$, 
% and the vector $f(t)\in\mbR^{N_V}$ represents external forces.
% The vector $y_{M} (t) \in \mbR^{N_p}$ represents boundary conditions and has non-zero elements only for those pressure volumes that are adjacent to a boundary and having a non-zero (possibly time-dependent) normal velocity component prescribed. $y_{M}(t)$ is obtained as a linear mapping $F_{M}$ from the actual boundary condition values $y_{bc}(t)$ as
% \begin{equation}
%     y_{M}(t) = F_{M}(y_{bc}(t)),\label{eq:mapping y_M y_bc}
% \end{equation}
% where $y_{bc}:[0,T]\rightarrow\mbR^{N_{V_{bc}}+N_{p_{bc}}}$ is the vector containing the velocity and pressure values at all boundary faces, as determined by the boundary conditions for $u_D$ and $p_{\infty}$. For details on $y_{bc}$ and $F_{M}$, we refer to \ref{app:discretiz of bc}.
The 
% In the
% right-hand side of the momentum equation 
% \eqref{eq:spacedisc navier stokes momentum}
% % is written as
% % can be decomposed as
% we can separate the pressure-dependent term
% 	\begin{align}
% 	F_h(V_h,p_h,\ybc,f) &= F_h^{CD}(V_h,\ybc,f) - G_h p_h
% 	\label{eq:mom rhs coarse decomp}
% % 	F_h^{CD}(V_h) &= -\tilde C_h(V_h)V_h + \nu D_hV_h + f_h.
% \end{align}
% Here,
% the 
term 
% $F_h^{CD}(V_h,\ybc,f)$
$F_h^{CD}(V_h(t),\ybc(t))$
contains the discretized convection, diffusion and forcing terms as well as all boundary condition contributions. The discretized forcing term $f(t)\in\mbR^{N_V}$ is assumed fixed and time-independent, hence neither $f(t)$ nor $t$ appear as an argument of $F_h^{CD}$.
% , which includes terms that result from time-dependent boundary conditions.
The term
$G_h p_h$ is the discretization of the (integrated) pressure gradient $\nabla p$ in 
\eqref{eq:conti mom cons} and $G_h\in\mbR^{N_V\times N_p}$ 
satisfies the discrete compatibility relation between the divergence and gradient operator
\begin{align}
	G_h = -M_h^T,
	\label{eq:spacedisc gradient divergence}
\end{align}
for all types of boundary conditions. This is a crucial property of the staggered grid discretization that will be used at several points in the construction of our new 
% mass-conserving 
energy-consistent
ROM approach for time-dependent boundary conditions.

To solve the system \eqref{eq:spacedisc navier stokes mass}-\eqref{eq:spacedisc navier stokes momentum}, we apply the discrete divergence operator $M_h$ to the $\Ohi$-premultiplied momentum equation \eqref{eq:spacedisc navier stokes momentum} and insert the mass equation \eqref{eq:spacedisc navier stokes mass} to get the Poisson equation for the pressure,
\begin{align}
    L_h p_h(t) := M_h\Ohi G_h p_h(t) = - \ddt y_M(t) + M_h\Ohi F_h^{CD}(V_h(t),\ybc(t)) .
	\label{eq:fom pressure possion eq}
\end{align}
As a result, 
% the pressure $p_h(t)$ in the momentum equation \eqref{eq:spacedisc navier stokes momentum} ensures that the velocity $V_h(t)$ satisfies the mass equation \eqref{eq:spacedisc navier stokes mass}.
when integrating the velocity $V_h(t)$ according to the momentum equation \eqref{eq:spacedisc navier stokes momentum}, the pressure $p_h(t)$ computed from \eqref{eq:fom pressure possion eq} ensures the velocity to satisfy the mass equation \eqref{eq:spacedisc navier stokes mass}.

Note that this approach requires the invertibility of $L_h$. As described in \cite{sanderse2012accuracy}, $L_h$ is singular for certain boundary conditions and the pressure $p_h$ is only defined up to an additional constant.
% but can be made invertible by imposing an additional constraint.
% but 
However, by adding a constraint, we can construct an invertible matrix $\bL\in\mbR^{N_p\times N_p}$ such that the unique solution to the regularized Poisson equation
\begin{align}
    \bL p_h = -\ddt y_M(t) + M_h\Ohi F_h^{CD}(V_h(t),\ybc(t))
    \label{eq:reg fom pressure poisson eq}
\end{align}
is also a solution to the Poisson equation \eqref{eq:fom pressure possion eq}. For boundary conditions for which $L_h$ is singular, we therefore compute $p_h$ via \eqref{eq:reg fom pressure poisson eq} to simulate the FOM \eqref{eq:spacedisc navier stokes mass}-\eqref{eq:spacedisc navier stokes momentum}.

To avoid the need to distinguish between cases with $L_h$ singular and invertible in our notation, 
% we will resort to the inverse $\Li$, 
we consider \eqref{eq:reg fom pressure poisson eq} the equation to be solved within the FOM simulation,
where we set $\bL := L_h$ for boundary conditions for which $L_h$ is invertible.
In both cases, the inverse of $\bL$ satisfies
\begin{align}
    L_h \Li y = y \h \text{ for all } y \in \Ima(L_h)
    % = \Im(M_h) \superset \Im(F_M).
    , 
    \h \text{ where } \h \Ima(F_M) \subset \Ima(M_h) = \Ima(L_h),
    \label{eq:li inverse property}
\end{align}
and 
% consequently
\begin{align}
    y^T \Lit z 
    % = y^T \Lit L_h \Li z = (L_h \Li y)^T \Li z 
    = y^T \Li z.
    \label{eq:li symmetry property}
\end{align}
In other words, within 
% $\Ima(M_h)$, 
the image space of $M_h$,
$\Li$ is symmetric and acts as an inverse to $L_h$.

% The term $F_h^{CD}(V_h,\ybc,t)$ contains the discretized convection, diffusion and forcing terms, which includes terms that result from time-dependent boundary conditions.

Suitable time discretizations of 
% the equations \eqref{eq:spacedisc navier stokes mass}-\eqref{eq:spacedisc navier stokes momentum} 
% this system
the system consisting of \eqref{eq:spacedisc navier stokes momentum} and \eqref{eq:reg fom pressure poisson eq}
are, for example, described in \cite{sanderse2012accuracy,sanderse2013energy}. 

% \section{Velocity-only ROM}
\section{Novel 
% mass-conserving, 
velocity-only and energy-consistent ROM for time-dependent boundary conditions}
\label{sec:velo only rom}

% \subsection{Helmholtz-Hodge-inspired decomposition}
\subsection{Lifting function inspired by Helmholtz-Hodge decomposition}
\label{sec:hhd inspired decomposition}
To construct a ROM for the FOM described in Section \ref{sec: a struc-pres fv disc}, we generalize the ideas in \cite[Section 4]{sanderse2020non} to time-dependent boundary conditions.
We decompose the FOM velocity,
\begin{align}
    V_h(t) = \Vhhom(t) + \Vhinhom(t),
    \label{eq:fom velocity decomp}
\end{align}
such that the time-dependent lifting function $\Vinhom(t)$ satisfies the mass equation,
\begin{align}
    M_h \Vhinhom(t) = y_M(t),
    \label{eq:inhom div free vhinhom}
\end{align}
and consequently, the homogeneous component $\Vhhom(t)$ satisfies the homogeneous mass equation,
% $M_h \Vhhom(t) = 0$. 
\begin{align}
    M_h \Vhhom(t) = 0.
    \label{eq:hom div free vhhom}
\end{align}

The choice of the lifting function $\Vinhom(t)$ is of pivotal importance. 
Inspired by the Helmholtz-Hodge decomposition \cite{bhatia2012helmholtz},
we choose the lifting function as
% We achieve \eqref{eq:inhom div free vhinhom}, by defining the lifting function as
% \begin{align}
%     \Vhinhom(t) = \Omega_h^{-1}G_hL_h^{-1}y_M(t)
%     = \Omega_h^{-1}G_hL_h^{-1}F_M\ybc(t).
%     \label{eq:def lifting function}
% \end{align}
\begin{align}
    \Vhinhom(t) = \Omega_h^{-1}G_h\Li y_M(t)
    = \Omega_h^{-1}G_h\Li F_M\ybc(t).
    \label{eq:def lifting function}
\end{align}
This choice has three decisive properties.
First, the lifting function indeed satisfies the mass equation with the inhomogeneous right-hand side, \eqref{eq:inhom div free vhinhom}. 
This property implies the homogeneous mass equation \eqref{eq:hom div free vhhom} for the homogeneous component $\Vhhom(t)$. 
We use this homogeneous mass equation 
to construct a velocity-only ROM in the next section. 
%  in Section \ref{sec:velo only rom subsec}.

Second, the lifting function is a linear function of the boundary condition vector $\ybc(t)$. We use this property to compute the ROM efficiently for arbitrary time-dependent boundary conditions in Section \ref{sec:effic rhs}.

Third, as in the Helmholtz-Hodge decomposition, the components in \eqref{eq:fom velocity decomp} are orthogonal, namely with respect to the $\Oh$-inner product,
% \begin{align}
%     \Vhhom^T\Omega_h\Vhinhom = \Vhhom^T\Omega_h\Omega_h^{-1}G_hL_h^{-1}y_M(t) = (M_h\Vhhom)^TL_h^{-1}y_M(t) = 0.
%     \label{eq:ortho of components 1}
% \end{align}
\begin{align}
    \Vhhom^T\Omega_h\Vhinhom = \Vhhom^T\Omega_h\Omega_h^{-1}G_h\Li y_M(t) = (M_h\Vhhom)^T\Li y_M(t) = 0.
    \label{eq:ortho of components 1}
\end{align}
This property guarantees the energy consistency of the ROM as described in Section \ref{sec:kin energy evol}.

\subsection{Velocity-only ROM}
\label{sec:velo only rom subsec}

% and
% % We 
% highlight two properties of this definition. First,
% to compute the lifting function $\Vinhom(t)$ at any point in time, the only information we need is the boundary condition vector $\ybc(t)$ at this point in time. Second, the lifting function is a linear function of this boundary condition vector.
% Hence, to compute the lifting function $\Vinhom(t)$ at a given time step, the only information we need is the mass equation right-hand side $y_M(t)$ at this time step.
% % As 
% Moreover, as
% explained in \eqref{eq:yM linear function of F_M} in \ref{app:discretiz of bc}, the mass equation boundary contribution $y_M(t)$ can be expressed as a linear function of the boundary condition vector $\ybc(t)$. Hence, we can write the lifting function $\Vinhom(t)$ as a linear function of $\ybc(t)$,
% \begin{align}
%     \Vinhom(t) = \Omega_h^{-1}G_hL_h^{-1}F_M\ybc(t) =: \Finhom \ybc(t),
% \end{align}
% with $\Finhom\in\mbR^{N_V\times \Nbc}$.

To approximate the homogeneous velocity $\Vhhom(t)$, we construct a POD basis with respect to the $\Omega_h$-inner product computed from snapshots $\Vhhom^j = V_h^j - \Vhinhom(t^j)$, where  $V_h^j$ is the time-discrete approximation of the velocity $V_h(t)$ defined by the FOM \eqref{eq:spacedisc navier stokes mass}-\eqref{eq:spacedisc navier stokes momentum} at time step $t^j$.
% As a result, 
Due to the homogeneous mass equation \eqref{eq:hom div free vhhom},
the resulting $\Oh$-orthonormal POD basis $\phihom\in\mbR^{N_V\times \Rhom}$ 
% with respect to the $\Omega_h$-inner product computed from snapshots $\Vhhom^n = V_h^n - \Vhinhom(t^n)$, with $V_h^n$ obtained from \eqref{eq:ex euler mass eq}-\eqref{eq:ex euler momentum eq}, 
satisfies
\begin{align}
    M_h\phihom = 0,
    \label{eq:m h phi hom = 0}
\end{align}
as described in \cite{sanderse2020non}.
Hence, we can construct the low-dimensional approximation $\Vrhom(t) = \phihom \ahom(t) \; \approx \Vhhom(t)$ such that the ROM $V_r(t) = \Vrhom(t) + \Vhinhom(t)$ satisfies the mass equation, 
% $M_hV_r(t) = y_M(t)$.
\begin{align}
    M_h V_r(t) = y_M(t).
    \label{eq:rom exact mass eq}
\end{align}

To determine the coefficient vector $\ahom(t)\in\mbR^{\Rhom}$, we replace $V_h(t)$ by $V_r(t)$ in the momentum equation \eqref{eq:spacedisc navier stokes momentum} and project onto $\phihom$,
\begin{align}
    \phihom^T \Omega_h\ddt V_r(t) = \ddt \ahom(t) + \phihom^T\Oh \ddt \Vinhom(t) =  \phihom^T F_h^{CD}(V_r(t),\ybc(t)) - \phihom^T G_hp_h(t).
    \label{eq:provisional rom ode}
\end{align}
Due to the divergence-gradient duality \eqref{eq:spacedisc gradient divergence} and property \eqref{eq:m h phi hom = 0}, the pressure term vanishes, so the ODE \eqref{eq:provisional rom ode} simplifies to
\begin{align}
    \ddt \ahom(t) + \phihom^T\Oh \ddt \Vinhom(t) = \phihom^T F_h^{CD}(\phihom \ahom(t) + \Vhinhom(t),\ybc(t))
    \label{eq:provisional rom ode2}.
\end{align}
This ODE does not contain the pressure $p_h$ anymore, so our ROM is velocity-only, i.e., the pressure is not needed to simulate the velocity. This property is very useful for two reasons. First, because we avoid additional computational costs. Second, because we avoid pressure-related $\inf$-$\sup$-instabilities \cite{ballarin2015supremizer}.

The steps taken so far also hold for other choices of lifting functions, e.g., \cite{mohebujjaman2017energy,grassle2019pod}. For the next simplification of ODE \eqref{eq:provisional rom ode2}, however, we exploit the specific definition of our lifting function, \eqref{eq:def lifting function},
namely the orthogonality property \eqref{eq:ortho of components 1}.
% This definition is 
% inspired by
% the Helmholtz-Hodge decomposition (HHD) \cite{bhatia2012helmholtz},
% and
% like the HHD, the components of our decomposition \eqref{eq:fom velocity decomp} are orthogonal, namely with respect to the $\Omega_h$-inner product,
% \begin{align}
%     \Vhhom^T\Omega_h\Vhinhom = \Vhhom^T\Omega_h\Omega_h^{-1}G_hL_h^{-1}y_M(t) = (M_h\Vhhom)^TL_h^{-1}y_M(t) = 0.
% \end{align}
Due to property \eqref{eq:m h phi hom = 0}, this orthogonality property is inherited by the POD basis $\phihom$, i.e.,
% $\phihom^T\Omega_h\Vhinhom(t) = 0$. 
\begin{align}
    \phihom^T\Omega_h\Vhinhom(t) = 0.
    \label{eq:ortho property}
    % \label{eq:ortho of components 2}
\end{align}
Due to this orthogonality, the left-hand side of \eqref{eq:provisional rom ode2} simplifies to $\ddt\ahom$, so the ODE for the ROM coefficient $\ahom(t)$ becomes
\begin{align}
    \ddt \ahom(t) = \phihom^T F_h^{CD}\left(\phihom\ahom(t) + 
    % \Oh^{-1}G_hL_h^{-1}F_M\ybc(t)
    \Vinhom(t)
    ,\ybc(t)\right)
    \label{eq:provisional rom ode3}
\end{align}
The initial condition is given by the projection of the FOM initial condition, $\ahom(0) = \phihom^TV_h(0)$.

% This property is pivotal for the derivation of the ROM kinetic energy evolution in Section \ref{sec:kin energy evol}.
% Moreover, the left-hand side of the ODE \eqref{eq:provisional rom ode2} simplifies to 
% % $\phihom^T\Omega_h \ddt\Vrhom(t) = $
% $\ddt\ahom$.

% \subsection{Efficient numerical simulation}
\subsection{Efficient evaluation of the ODE right-hand side}
\label{sec:effic rhs}

To simulate our novel ROM, 
% numerically, 
we have to discretize the ODE
\eqref{eq:provisional rom ode3}
% \begin{align}
%     \ddt \ahom(t) = \phihom^T F_h^{CD}\left(\phihom\ahom(t) + 
%     % \Oh^{-1}G_hL_h^{-1}F_M\ybc(t)
%     \Vinhom(t)
%     ,\ybc(t)\right)
%     \label{eq:provisional rom ode3}
% \end{align}
in time and evaluate the right-hand side. As a naive approach, we could first evaluate the 
% term 
FOM-like right-hand side
$F_h^{CD} (\phihom \ahom(t) +\Vinhom(t),\ybc(t))$ and then perform the projection onto $\phihom$, both during the online phase of the ROM simulation. However, the computational costs would
% scale at least linearly in $N_V$ 
depend on the FOM dimensions
and hence deteriorate the computational complexity of the ROM in the online phase substantially.

Common ways of avoiding such prohibitive costs are
hyper-reduction techniques such as the discrete interpolation method (DEIM) \cite{chaturantabut2009discrete} and gappy POD \cite{brunton2019data}. 
% These methods approximate the evaluation of a nonlinear function by evaluating this nonlinear function only in a small number of entries. 
Instead of evaluating the complete FOM-dimensional vector $F_h^{CD} (\phihom \ahom(t) +\Vinhom(t),\ybc(t))$, these methods only require a small number of entries to construct an approximation of this vector. %Consequently, these methods improve the efficiency if computing only few entries is less expensive than computing the complete vector.
% If the evaluation of only a small number of entries is significantly less expensive than the evaluation of the complete FOM-dimensional vector, we
However, the definition of our lifting function $\Vinhom(t)$ \eqref{eq:def lifting function} contains the inverse 
% $L_h^{-1}$,
$\Li$,
requiring the solution of a
% (regularized) 
Poisson equation with computational complexity $\bigo{N_p^2}$. % so that evaluating a few entries of $\Vinhom$ is not straightforward. 
%Computing this inverse explicitly is typically prohibitively time- and memory-consuming, so, in practice, matrix multiplications $L_h^{-1}b$ with $b\in\mbR^{N_p}$, are computed instead by solving the Poisson equation $L_h x = b$.
Consequently, sampling a few entries of $\Vinhom(t)$ is already prohibitively expensive.

Instead of applying hyper-reduction to the right-hand side of \eqref{eq:provisional rom ode3} as a whole, we therefore propose to approximate only the boundary condition vector $\ybc(t)$.
In detail,
% To solve these issues, 
we propose to replace the boundary condition vector $\ybc(t)$ in the ODE \eqref{eq:provisional rom ode3} and in the definition of the lifting function \eqref{eq:def lifting function} with a linear approximation 
% $\tyb(t) = \phibc \abc(t) \; \approx \ybc(t)$ 
\begin{align}
    \tybc(t) = \phibc \abc(t) \; \approx \ybc(t),
    \label{eq:boundary condition approximation}
\end{align}
with $\phibc\in\mbR^{\Nbc\times \Rbc}, \abc(t)\in\mbR^\Rbc$ and $\Rbc \ll \Nbc$. 
One option to construct this approximation are hyper-reduction techniques such as DEIM \cite{chaturantabut2009discrete}.
% This approximation could be constructed via hyper-reduction techniques such as DEIM \cite{chaturantabut2009discrete}. 
% For this purpose, 
In that case,
the basis $\phibc$ would be computed via a POD of snapshots of $\ybc(t)$ during the offline phase. The coefficient vector $\abc(t)$ would be computed based on the evaluation of only $\Rbc$ entries of $\ybc(t)$ with complexity $\bigo{\RbcII}$ at each time step during the online phase.
However, such a DEIM coefficient vector is only an approximation to the optimal coefficient. 
% In fact, 
Actually,
the best approximation to $\ybc(t)$ within the space spanned by the basis $\phibc$ is given by $\phibc\phibc^T\ybc(t)$. 
Hence, we can improve the accuracy of the approximation $\tybc(t)$ 
by computing the coefficient vector via the projection onto $\phibc$, $\abc(t) = \phibc^T \ybc(t)$, instead of performing DEIM. 
This computation would be prohibitively expensive when performed during the online phase.
% While this computation would be prohibitively expensive when performed during the online phase, we can 
% use 
% \textit{eager precomputation} 
However, we can maintain the online efficiency by
computing $\abc(t^j)$ for all required time steps $t^j$ already in the offline phase.
This approach requires that $\ybc(t)$ is known a priori, which we assume in this work. If $\ybc(t)$ would not be known a priori, we suggest to resort to computing $\abc(t)$ via DEIM.
% to move the computational cost to the offline phase. Storing the coefficient vector $\abc(t)\in\mbR^\Rbc$ for all time steps requires a memory capacity of only $\bigo{\Rbc\Nttotal}$ which can be further reduced if we precompute only a smaller number of time steps and interpolate during the online phase. Such an interpolative approach also mitigates the need to know all time steps a priori.

The replacement of the boundary condition vector $\ybc(t)$ by the approximation \eqref{eq:boundary condition approximation} introduces an error in the ROM.
% By replacing the boundary condition vector $\ybc(t)$ by the approximation \eqref{eq:boundary condition approximation},
% we introduce an error to the ROM.
% This approximation introduces an additional error in the ROM. 
In particular, 
% In detail,
the resulting ROM 
$V_r(t) = \phihom \ahom(t) + \tVinhom(t)$ 
% $V_r$
does not satisfy the mass equation with respect to the original boundary conditions \eqref{eq:rom exact mass eq}, but with respect to the approximated 
% boundary 
boundary conditions,
% right-hand side,
\begin{align}
    M_h V_r(t) = \tilde y_M := F_M \tybc(t).
\end{align}
% and the coefficient $\tahom(t)$ 
The coefficient $\ahom(t)$ is defined by the ODE
\begin{align}
    \ddt \ahom(t) = \phihom^T F_h^{CD}\left(\phihom\ahom(t) + 
    % \Oh^{-1}G_hL_h^{-1}F_M\ybc(t)
    % \Finhom\phibc\abc(t)
    \tVinhom(t)
    ,\phibc\abc(t)\right)
    % .
    ,
    \label{eq:approximated rom ode}
\end{align}
where the approximated lifting function $\tVinhom(t)$ is defined by
% \begin{align}
%     \tVinhom(t) = \Oh^{-1} G_hL_h^{-1} F_M\phibc\abc(t) 
%     =: \Finhom\abc(t).
%     \label{eq:def approx lifting function}
% \end{align}
\begin{align}
    \tVinhom(t) = \Ohi G_h\Li F_M\phibc\abc(t) 
    =: \Finhom\abc(t).
    \label{eq:def approx lifting function}
\end{align}
The main advantage of this approximation is 
% that 
% this approximation
% reduces 
the drastic reduction of
the number of (regularized) Poisson solves required for computing
% $\Vrinhom(t)$ 
the lifting function
% drastically 
due to the linear dependence on $\ybc(t)$ in the definition \eqref{eq:def lifting function}. 
% For the computation of the approximated lifting function
% \begin{align}
%     \tVinhom(t) = \Oh^{-1} G_hL_h^{-1} F_M\phibc\abc(t) 
%     =: \Finhom\phibc\abc(t), 
% \end{align}
Indeed,
to precompute the matrix $\Finhom\in\mbR^{N_V\times\Rbc}$, we
need to solve only $\Rbc$ (regularized) Poisson equations for the $\Rbc$ columns of the matrix $F_M\phibc$. 

Note that this approximated lifting function satisfies the same orthogonality property as the exact lifting function, \eqref{eq:ortho property},
\begin{align}
    \phihom^T \Oh\tVinhom(t) = 0
    \label{eq:ortho of components 3}.
\end{align}
We will use this orthogonality property in Section \ref{sec:kin energy evol} to derive energy consistency.

\subsection{Exact offline decomposition using matricized third-order tensors}
To evaluate the right-hand side of the approximated ODE \eqref{eq:approximated rom ode}
efficiently, we extend the offline precomputation approach in \cite{sanderse2020non}.
This approach
exploits the low-dimensional structure of the terms $\phihom\ahom(t)$, $\Finhom\abc(t)$ and $\phibc\abc(t)$ as well as the fact that 
% $F_h^{CD}$ consists of polynomials of its arguments not higher than second order. 
the nonlinearities in $F_h^{CD}$ are only second order polynomials.
In detail,
we formulate these second order polynomials as precomputable ``matricized" third order tensors. For example, we can write the term quadratic in $\phihom\ahom(t)$ as $\bar C_2 \, (\phihom\ahom(t)\otimes \phihom\ahom(t))$ with $\bar C_2\in\mbR^{N_V\times N_V^2}$ and the Kronecker product $\otimes$. Then, we find
\begin{align}
    \phihom^T \bar C_2 \, (\phihom\ahom(t)\otimes \phihom\ahom(t)) = \phihom^T \bar C_2 \, (\phihom\otimes\phihom)\, (\ahom(t)\otimes\ahom(t))
\end{align}
where $\phihom^T \bar C_2 \, (\phihom\otimes\phihom)\in\mbR^{\Rhom\times\RhomII}$ can be precomputed offline. Similarly, we can precompute operators in $\mbR^{\Rhom\times \RhomII}$ and $\mbR^{\Rhom\times(\Rhom\Rbc)}$ for the term quadratic in $\phibc\abc(t)$ and the mixed term linear in both, $\phihom\ahom(t)$ and $\phibc\abc(t)$, respectively.
Hence, the resulting method has an online complexity scaling in $\bigo{\RhomIII + \RhomII\Rbc + \Rhom\RbcII}$.
% As a more efficient alternative
% If larger numbers of modes are required that prohibit this cubic scaling, we suggest to incorporate energy-conserving hyper-reduction

% The contribution of the force term $f(t)$ can be dealt with either via hyper-reduction or \textit{eager precomputation}.

\subsection{Pressure recovery}
\label{sec:pres recovery}

As discussed in Section \ref{sec:velo only rom subsec}, our proposed ROM is velocity-only, so no computation of the pressure is needed to simulate the velocity. However, we can compute the pressure $p_r(t)$ of the ROM as a post-processing step by inserting the ROM velocity $V_r(t)$ and the boundary condition approximation $\ybc(t)$ into the (regularized) FOM Poisson equation 
% \eqref{eq:fom pressure possion eq},
\eqref{eq:reg fom pressure poisson eq}
% \begin{align}
%     L_h p_r(t) = - \ddt \tilde y_M(t) + M_h \Ohi F_h^{CD}(V_r(t),\tybc(t)).
%     \label{eq:vo rom poisson eq}
% \end{align}
\begin{align}
    \bL p_r(t) = - \ddt \tilde y_M(t) + M_h \Ohi F_h^{CD}(V_r(t),\tybc(t)).
    \label{eq:vo rom poisson eq}
\end{align}
This pressure is consistent in the sense that 
the (regularized) Poisson equation \eqref{eq:vo rom poisson eq} has the same structure as the (regularized) FOM Poisson equation \eqref{eq:fom pressure possion eq}. Consequently,
the ROM pressure $p_r(t)$ converges to the FOM pressure $p_h(t)$ if the ROM velocity $V_r(t)$ converges to the FOM velocity $V_h(t)$ and the boundary condition approximation $\tybc(t)$ converges to the exact boundary condition vector $\ybc(t)$.

\section{Kinetic energy evolution}
    \label{sec:kin energy evol}
    
A distinctive property of 
our proposed ROM 
% as given in \eqref{eq:approximated rom ode},
is its kinetic energy consistency, i.e., the ROM kinetic energy evolution has the same structure as the FOM kinetic energy evolution. We define the FOM kinetic energy as $K_h(t) = \frac{1}{2}\|V_h(t)\|_{\Omega_h}^2$, leading to the kinetic energy 
% balance
evolution
\begin{align}
    \ddt K_h(t) = 
    % V_h^T\Omega_h\ddt V_h = V_h^TF_h^{CD}(V_h,\ybc) - V_h^TG_h p_h =
    V_h(t)^TF_h^{CD}(V_h(t),\ybc(t)) + y_M(t)^Tp_h(t),
    \label{eq:fom kin energy balance}
\end{align}
% \begin{align}
%     \frac{K_r^{n+1}-K_r^n}{\Delta t} = V_h^n^T\Oh \frac{V_h^{n+1}-V_h^n}{\Delta t}= V_h^n^T F_h^{CD}(V_h^n,\ybc(t^n),f) 
% \end{align}
% time discrete energy evolution only works out for energy-conserving (=> implicit) time discretization schemes
using the momentum equation \eqref{eq:spacedisc navier stokes momentum}
% , the RHS decomposition \eqref{eq:mom rhs coarse decomp}
and the mass equation \eqref{eq:spacedisc navier stokes mass}. 

The key property leading to a kinetic energy evolution with the same structure for our proposed ROM is the orthogonality of $\Vrhom(t)$ and $\tVinhom(t)$, given by equation \eqref{eq:ortho of components 3}.
Due to 
% the orthogonality of $\Vrhom$ and $\tVinhom$ \eqref{eq:ortho of components 3},
this orthogonality
the 
% velocity-only 
ROM kinetic energy can be decomposed into two terms,
\begin{align}
    K_r(t) &= \frac{1}{2} \|V_r(t)\|_\Oh^2 
    % = \frac 1 2 \left(\phihom\ahom(t) +\Vinhom(t)\right)^T\Oh\left(\phihom\ahom(t) +\Vinhom(t)\right) \\
    % &= \frac 1 2 \left[\ahom(t)^T\phihom^T\Oh\phihom\ahom(t) + \Vinhom(t)^T\Oh\Vinhom(t)\right] 
    = \frac 1 2 \|\ahom(t)\|_2^2 + \frac 1 2 \|\tVinhom(t)\|_\Oh^2,
    \label{eq:vo rom energy decomp}
\end{align}
%where we exploited the orthogonality property \eqref{eq:ortho property} of the HHD-inspired decomposition and the $\Oh$-orthonormality of $\phihom$. 
The time-derivative of the first term (the energy of the homogeneous component) is given by
\begin{align}
    % \ddt \frac 1 2 \|\ahom\|_2^2 = \ahom^T\ddt \ahom = \ahom^T\phihom^TF_h^{CD}(V_r,\tybc),
    \ddt \frac 1 2 \|\ahom(t)\|_2^2 = \ahom(t)^T\ddt \ahom(t) = \ahom(t)^T\phihom^TF_h^{CD}(V_r(t),\tybc(t)),
\end{align}
using the momentum equation \eqref{eq:approximated rom ode}.
Furthermore, the time-derivative of the second term (the energy of the approximated lifting function \eqref{eq:def lifting function}) 
% is given by
can be written as
% \begin{align}
%     \ddt \frac 1 2 \|\tVinhom\|_\Oh^2 
%     &= \tVinhom^T\Oh \ddt \tVinhom = \left(\Ohi G_h L_h^{-1}\tilde{y}_M\right)^T\Oh \ddt \left(\Ohi G_h L_h^{-1}\tilde{y}_M\right) \\
%     &= -\tilde y_M^T L_h^{-1} M_h\Ohi G_h L_h^{-1}\ddt \tilde{y}_M = -\tilde y_M^T L_h^{-1} \ddt \tilde{y}_M, 
% \end{align}
\begin{align}
    \ddt \frac 1 2 \|\tVinhom(t)\|_\Oh^2 
    &= \tVinhom(t)^T\Oh \ddt \tVinhom(t) = \left(\Ohi G_h \Li\tilde{y}_M(t)\right)^T\Oh \ddt \left(\Ohi G_h \Li\tilde{y}_M(t)\right) \\
    &= -\tilde y_M(t)^T \bL^{-T} M_h\Ohi G_h \Li\ddt \tilde{y}_M(t) 
    % = - (L_h \Li \tilde y_M)^T \Li \ddt \tilde y_M
    = -\tilde y_M(t)^T \Li \ddt \tilde{y}_M(t), 
\end{align}
    % \intertext{
    where we use the divergence-gradient duality \eqref{eq:spacedisc gradient divergence}, the symmetry of $L_h$ and $\Oh$, the definition of $L_h$ \eqref{eq:fom pressure possion eq}, and the inverse and symmetry properties of $\Li$ \eqref{eq:li inverse property}-\eqref{eq:li symmetry property}.
    
    Insertion of the 
    % (regularized) 
    Poisson equation for the ROM pressure \eqref{eq:vo rom poisson eq} leads to
    % \begin{align}
    %     \ddt \frac 1 2 \left\|\tVinhom\right\|_\Oh^2 
    %     &= -\tilde y_M ^T \left[ - p_r + L_h^{-1}M_h \Ohi F_h^{CD}(V_r,\tybc) \right] \\
    %     &= \tilde y_M^T p_r - \tVinhom^T F_h^{CD}(V_r,\tybc).
    % \end{align}
    \begin{align}
        \ddt \frac 1 2 \left\|\tVinhom(t)\right\|_\Oh^2 
        &= -\tilde y_M(t) ^T \left[ - p_r(t) + \Li M_h \Ohi F_h^{CD}(V_r(t),\tybc(t)) \right] \\
        &= \tilde y_M(t)^T p_r(t) - \tVinhom(t)^T F_h^{CD}(V_r(t),\tybc(t)).
    \end{align}
    % We add and substract the term $\tilde{y}_M^T L_h^{-1}\left[M_h\Ohi F_h^{CD}(V_r,\tybc) \right]$, so that we can insert 
    % % the time-continuous counterpart of 
    % the Poisson equation for the pressure \eqref{eq:fom pressure possion eq} given the velocity $V_r$ and the boundary condition vector $\tybc$,
    % }
% \begin{align}
%     \ddt \frac 1 2 \|\tVinhom\|_\Oh^2
%     &= \tilde{y}_M^T L_h^{-1} \left[ -\ddt \tilde y_M + M_h\Ohi F_h^{CD}(V_r,\tybc) - M_h\Ohi F_h^{CD}(V_r,\tybc)\right] \\
%     &= \tilde y_M^T p_h 
%     -\tilde y_M^T L_h^{-1} M_h\Ohi F_h^{CD}(V_r,\tybc).
%     % +\left(\Oh^{-1}G_hL^{-1}\tilde{y}_M\right)^T F_h^{CD}(V_r,\tybc,f).
% \end{align}
    % \intertext{
    Here, we have again used
    % Again using 
    the divergence-gradient duality \eqref{eq:spacedisc gradient divergence}, the symmetry of $L_h$ and $\Oh$, the symmetry property of $\Li$ \eqref{eq:li symmetry property}, 
    % we identify 
    and identified
    the definition of the lifting function \eqref{eq:def lifting function},
    \begin{align}
        -\tilde y_M(t)^T \Li M_h\Ohi 
        = -\tilde y_M(t)^T \Lit M_h\Ohi
        = \left(\Ohi G_h \Li\tilde{y}_M(t) \right)^T = \tVinhom(t)^T.
    \end{align}
%     so we find
%     % . Applying the definition of the lifting function \eqref{eq:def lifting function} one more time, we find
%     % }
%     \begin{align}
%     \ddt \frac 1 2 \|\tVinhom\|_\Oh^2 
%     &= \tilde{y}_M^Tp_h + \Vinhom^T F_h^{CD}(V_r,\tybc,f).
% \end{align}
% Hence, 
Altogether,
the kinetic energy evolution of the ROM can be written as
\begin{align}
    \ddt K_r(t) 
    % = [\phihom\ahom +\tVinhom]^T F_h^{CD}(V_r,\tybc,f) + \tilde{y}_M^Tp_h 
    = V_r(t)^T F_h^{CD}(V_r(t),\tybc(t)) + \tilde{y}_M(t)^T p_r(t).
    \label{eq:rom kin energy balance}
\end{align}
This
ROM energy evolution has exactly the same structure as the FOM energy evolution \eqref{eq:fom kin energy balance}, which is our definition of consistency.

% In \cite{sanderse2020non}, 
% such a structural correspondence of FOM and ROM kinetic energy evolution for homogeneous boundary conditions is used to derive nonlinear stability of the ROM: both the FOM and ROM energy are bounded.
% In contrast, we cannot derive such bounds for arbitrary boundary conditions. However, for specific boundary conditions
% such bounds can be derived, e.g. \cite[Book II, Chapter 11]{foias1978remarques} (on the continuous level). Thanks to the structural correspondence of equations \eqref{eq:fom kin energy balance} and \eqref{eq:rom kin energy balance}, a bound on the FOM energy implies a bound on the ROM energy. Hence, if the FOM energy is bounded for the boundary conditions $\tybc$, then the ROM energy is bounded, too, and the ROM is nonlinearly stable.

% As stressed in Section \ref{sec:velo only rom}, the orthogonality property \eqref{eq:ortho property} (see also \eqref{eq:ortho of components 1} and \eqref{eq:ortho of components 3})
% % of the HHD-inspired decomposition 
% is pivotal for the energy consistency described here. 
As indicated above, the orthogonality property \eqref{eq:ortho of components 3} is key to achieve this energy consistency.
We exploited this property in \eqref{eq:vo rom energy decomp} to eliminate
the term $\tVinhom(t)^T \Oh \phihom \ahom(t)$.
If property \eqref{eq:ortho property} would not hold, 
% these terms does in general not vanish. 
this term would not vanish and the energy evolution would contain the terms $\tVinhom(t)^T\Oh\phihom\ddt \ahom(t)$ and $(\phihom\ahom(t))^T\Oh \ddt \tVinhom(t)$.
This would be problematic because they
% This term 
lack clear physical interpretation, create a mismatch between ROM and FOM energy evolution and hence spoil the energy consistency.

% To obtain a energy evolution corresponding to \eqref{eq:rom kin energy balance} on the time discrete level, we would require energy-conserving Runge-Kutta methods. Since these methods are implicit, they can be prohibitively expensive. However, numerical experiments in \cite{sanderse2020non} showed that high-order explicit Runge-Kutta methods can have negligible energy errors.

On the time-discrete level, an energy evolution corresponding to \eqref{eq:rom kin energy balance} can be, for example, obtained using energy-conserving Runge-Kutta methods \cite{sanderse2013energy}. Since these methods are implicit, they can be prohibitively expensive. However, numerical experiments in \cite{sanderse2020non} have shown that high-order explicit Runge-Kutta methods can have negligible energy errors. If the time step size is not dominated by stability but by accuracy such explicit schemes can be a useful alternative to exactly energy-conserving implicit schemes.
\section{Equivalent velocity-pressure ROM}
\label{sec:equivalent velo pres rom}

\subsection{Velocity-pressure ROM}
\label{sec:vp rom}

Our proposed ROM has the appealing property that it is velocity-only. However, we can interpret this ROM as a velocity-pressure ROM, as well. This interpretation provides further insight in the relation between time-dependent boundary conditions and the role of the pressure.

% Although our proposed ROM is velocity-only, we have discovered this ROM to be equivalent to a velocity-pressure ROM. This equivalence provides crucial insights for understanding and comparing different ROMs.

The velocity-pressure ROM that we consider
comprises the velocity approximation $
%\Vrvp(t) = 
V_r(t) = \Phi a(t) \;\approx V_h(t)$ with the $\Oh$-orthonormal basis $\Phi\in\mbR^{N_V\times R_V}$ and the coefficient vector $a(t)\in\mbR^{R_V}$, and the pressure approximation $p_r(t) = \Psi b(t) \;\approx p_h(t)$ with the 
% orthonormal 
basis $\Psi\in\mbR^{N_p\times R_p}$ and the coefficient vector $b(t)\in\mbR^{R_p}$. 
In the next section, we will specify the choice of bases such that the resulting velocity-pressure ROM is equivalent to the velocity-only ROM proposed in Section \ref{sec:velo only rom}. In general, however, the bases can also be chosen differently.
% We will specify the choice of bases in the following sections, but note that $\Phi$ is generally different from $\phihom$ in Section \ref{sec:hhd inspired decomposition}. 

We insert these approximations in the mass equation \eqref{eq:spacedisc navier stokes mass} projected by the pressure basis $\Psi$ and in the momentum equation \eqref{eq:spacedisc navier stokes momentum} projected by the velocity basis $\Phi$, leading to:
\begin{align}
    % \hat M a(t) :=
    \Psi^TM_h \Phi a(t) &= \Psi^T 
    \tilde
    y_M(t) := \Psi^T F_M \tybc(t)
    % =: \hat y_M(t) 
    \label{eq:vp rom mass eq}
    \\
    \Phi^T\Oh \Phi \ddt a(t) = \ddt a(t) &= \Phi^T F_h^{CD}(\Phi a(t), 
    % \ybc(t)
    \tybc(t)
    ) - \Phi^T G_h \Psi b(t).
    \label{eq:vp rom momentum eq}
\end{align}
% The initial condition is defined via projection of the FOM initial condition, $a(0) = \Phi^T\Oh V_h(0)$.
The initial condition $a(t^0)=a^0$ is defined via the constrained optimization problem of minimizing $\|\Phi a^0 - V_h(t^0)\|_\Oh$ such that 
% $\Psi^T M \Phi a^0 = \Psi^T y_M(t^0)$.
$a^0$ satisfies the projected mass equation \eqref{eq:vp rom mass eq} at $t=t^0$.

Analogously to the FOM described in Section \ref{chap: intro to sanderse}, we solve this system by solving the reduced Poisson equation for $b(t)$,
\begin{align}
    L_r b(t) := \Psi^TM_h\Phi \Phi^T G_h \Psi b(t) = - \Psi^T \ddt 
    \tilde
    y_M(t) + \Psi^T M_h \Phi \Phi^T F_h^{CD}(\Phi a(t),
    % \ybc(t)
    \tybc(t)
    ).
    \label{eq:vp rom poisson eq}
\end{align}
Then, due to the pressure term in the ODE \eqref{eq:vp rom momentum eq}, $a(t)$ computed from this ODE is guaranteed to satisfy the projected mass equation \eqref{eq:vp rom mass eq}. In contrast to the FOM, the system matrix $L_r$ of this reduced Poisson equation is not high-dimensional and sparse, but low-dimensional and dense. 
% To solve this system, we precompute a LU-factorization
Hence, the system is solved (by employing a precomputed LU-factorization) 
% causing
with
a computational complexity of $\bigo{R_p^2}$ at each time step in the online phase.
This approach requires the invertibility of the system matrix $L_r = \Psi^TM_h\Phi \Phi^T G_h \Psi$. We will guarantee this invertibility by choosing the bases $\Phi$ and $\Psi$ suitably, as explained in the next section.

Analogously to the discussion in Section \ref{sec:effic rhs} for the velocity-only ROM, computing the right-hand side of the ODE \eqref{eq:vp rom momentum eq} for the exact boundary condition vector $\ybc(t)$ would in general be prohibitively expensive. 
% In contrast to that velocity-only ROM, the application of hyper-reduction methods to this velocity-pressure ROM would result in an efficient ROM since no lifting function spoils the efficiency gain. 
% However, with regard to the ROM equivalence that we show in the next section, 
Hence,
we have replaced the exact boundary condition vector $\ybc(t)$ in \eqref{eq:vp rom mass eq}-\eqref{eq:vp rom poisson eq} by the lower-dimensional approximation $\tybc(t) := \phibc \abc(t)$ as described in Section \ref{sec:effic rhs}.
% in 
% the velocity-pressure ROM
% \eqref{eq:vp rom mass eq}-\eqref{eq:vp rom poisson eq}.
% , instead of using hyper-reduction.

\subsection{Equivalence to the proposed velocity-only ROM}
\label{sec:eq to the vo rom}

% thoughts
% \begin{itemize}
%     \item Don't want to make unnecessarily general theorem: Don't want to introduce abstract variables, e.g., for $\phihom$ and $\phiinhom$.
%     \item[$\rightarrow$] idea: We need requirement that $\spannnn(\Phi) = \spannnn([\phihom \;\; \phiinhom]$ in any case, implying that there exists an orthogonal matrix $Q\in\mbR^{R_V\times R_V}$ such that $\Phi = [\phihom \;\; \phiinhom] Q$.
% \end{itemize}

% proof steps
% \begin{enumerate}
%     \item projected implies exact mass equation
%     \item exact mass equation has unique solution -> inhomogeneous components are equivalent
%     \item homogeneous components satisfy equivalent ODE
%     \item initial conditions are equivalent
% \end{enumerate}

The following theorem states
that the velocity-pressure ROM \eqref{eq:vp rom mass eq}-\eqref{eq:vp rom momentum eq} is equivalent to the velocity-only ROM \eqref{eq:approximated rom ode}-\eqref{eq:def approx lifting function} if we choose the bases $\Phi$ and $\Psi$ appropriately.

% how we have to choose the bases of the velocity-pressure ROM described in the previous section to get the same velocity as for our velocity-only ROM proposed in Section \ref{sec:velo only rom}.
% % Two theorems as on poster

\begin{mytheo}[\bf Equivalence of velocity-only and velocity-pressure ROM]
\label{eq:theo rom equivalence}
\ \\
If the $\Oh$-orthonormal velocity basis $\Phi\in\mbR^{N_V\times R_V}$ and the pressure basis $\Psi \in\mbR^{N_p\times R_p}$ are chosen such that
\begin{itemize}
% \begin{enumerate}
    \item[a)]
    the image space of $\Phi$ is equivalent to the union of the image spaces of
    % $\Phi$ spans the same subspace as the 
    % union of the spaces spanned by  
    % $[\phihom\;\;\Finhom]$  with 
    $\phihom$ and $\Finhom$ as defined in Section \ref{sec:velo only rom},
    \item[b)]
    for all $t\in[0,T]$, there exists a solution to the unprojected mass equation
    \begin{align}
        M_h \Phi a(t) = y_M(t) 
        % := F_M \ybc(t)
        ,
        \label{eq:vp rom mass eq exact}
    \end{align}
    \item[c)]
    $R_p = \rank(M_h\Phi)$,
    \item[d)]
     and
    $\Psi^T M_h \Phi$ has rank $R_p$,
    % the same rank as $M_h \Phi$,
    \end{itemize}
% \end{enumerate}
then the velocity of the velocity-pressure ROM, $\Vrvp(t) = \Phi a(t)$, is for all $t\in[0,T]$ equivalent to the velocity of the velocity-only ROM, $\Vrvo(t) = \phihom \ahom(t) + \tVinhom(t)$ described in Section \ref{sec:velo only rom}.
\end{mytheo}

\begin{proof}
    First of all, Assumption (a) implies that there exists an orthonormal matrix $Q\in\mbR^{R_V\times R_V}$ such that $\Phi = [\phihom \;\; \phiinhom]Q$. Here, $\phiinhom\in\mbR^{N_V\times \Rinhom}$ with $\Rinhom = \rank(\Finhom)$ is an $\Oh$-orthonormal basis with the same image space as $\Finhom$ which can be computed from $\Finhom$, e.g., via a Gram-Schmidt-orthonormalization. We can decompose the matrix $Q = \begin{bmatrix} Q_1\\ Q_2 \end{bmatrix}$ with $Q_1\in\mbR^{\Rhom \times R_V}$ and $Q_2 \in\mbR^{\Rinhom \times R_V}$. Then, we write the velocity-pressure ROM velocity as 
    \begin{align}
        \Vrvp(t) = \Phi a(t) = \phihom Q_1 a(t) + \phiinhom Q_2 a(t) =: \phihom a_1(t) + \phiinhom a_2(t).
        \label{eq:vp velocity decomp}
    \end{align}
    We proceed to show that the velocity components of the two ROMs are equivalent, namely
    \begin{align}
        a_1(t) = \ahom(t) \h \text{ for all } t\in[0,T],
        \label{eq:hom equivalence}
    \end{align}
    and
    \begin{align}
        \phiinhom a_2(t) = \tVinhom(t) \h \text{ for all } t\in[0,T].
        \label{eq:inhom equivalence}
    \end{align}
    
    For this purpose, we note that the decomposition \eqref{eq:vp velocity decomp} separates the modes in $\phihom$, which lie in the kernel of $M_h$, from those in $\phiinhom$, which lie outside the kernel of $M_h$. Consequently, equations \eqref{eq:vp rom mass eq} and \eqref{eq:vp rom mass eq exact} simplify to
    \begin{align}
        \Psi^T M_h \phiinhom a_2(t) &= \Psi^T \tilde y_M(t),
        \label{eq:simp vp rom mass eq}
    \intertext{ and }
        M_h \phiinhom a_2(t) &= \tilde y_M(t),
        \label{eq:simp vp rom mass eq exact}
    \end{align}
    respectively.
    Let, for a fixed $t$, $\hat a_2$ and $\bar a_2$ be the solutions to \eqref{eq:simp vp rom mass eq} and \eqref{eq:simp vp rom mass eq exact}, respectively (their existence follows from Assumption (b)). Then, premultiplying \eqref{eq:simp vp rom mass eq exact} by $\Psi^T$, we find
    \begin{align}
        \Psi^T M_h \phiinhom \bar a_2 = \Psi^T y_M(t) = \Psi^T M_h\phiinhom \hat a_2.
        \label{eq:proj equals exact}
    \end{align}
    By construction of $\phiinhom$, we have $\rank(M_h\Phi) = \rank(M_h\phiinhom) = \Rinhom$. Because of Assumption (c), $R_p = \Rinhom$, and because of Assumption (d), $\rank(\Psi^TM_h\phiinhom)=\rank(\Psi^TM_h\Phi) = \Rinhom$. Hence, $\Psi^TM_h\phiinhom\in\mbR^{\Rinhom\times \Rinhom}$ is invertible, so \eqref{eq:proj equals exact} implies $\bar a_2 = \hat a_2$. Hence, the solution $\hat a_2$ to the projected mass equation \eqref{eq:simp vp rom mass eq} also satisfies the  unprojected mass equation \eqref{eq:simp vp rom mass eq exact}. This solution $\hat a_2$ to the unprojected mass equation is unique. To show this uniqueness, let $\check a_2$ be a second solution to \eqref{eq:simp vp rom mass eq exact}. Then, we find
    \begin{align}
        M_h\phiinhom(\hat a_2 - \check a_2) = y_M(t) - y_M(t) = 0.
    \end{align}
    Since $M_h\phiinhom$ has full column rank, we infer $\hat a_2 - \check a_2 = 0$, so the two solutions are equivalent. This uniqueness implies the equivalence of the inhomogeneous velocity components \eqref{eq:inhom equivalence} because $\Finhom$ has the same image space as $\phiinhom$ and both terms satisfy the unprojected mass equation for all $t\in[0,T]$.

    What remains to be shown is the equivalence of the homogeneous velocity components \eqref{eq:hom equivalence}. For this purpose, we premultiply the momentum equation \eqref{eq:vp rom momentum eq} by $Q_1$ to find
    \begin{align}
        \ddt Q_1 a(t) = \ddt a_1(t) = \phihom^T F_h^{CD}(\phihom a_1(t)+\tVinhom(t),\tybc(t)),
        \label{eq:vp rom hom ode}
    \end{align}
    where $Q_1 \Phi^T = \phihom^T$ and the pressure term $\phihom^TG_h \Psi b(t)$ disappears due to equation \eqref{eq:m h phi hom = 0}. This ODE is equivalent to the velocity-only ODE \eqref{eq:approximated rom ode}. As both ODEs have the same initial condition, they also have the same solution, so we find \eqref{eq:hom equivalence}.

    Altogether, we hence find 
    % for all $t\in[0,T]$
    \begin{align}
        \Vrvp(t) 
        % = \phihom a_1(t) + \phiinhom a_2(t) = \phihom \ahom(t) + \tVinhom(t) 
        = \Vrvo(t)
        % \\ 
        \h \text{ for all } t\in[0,T]
        .
    \end{align}
    \qed
\end{proof}
The requirements of the theorem are met, for example, by the choice $\Phi = [\phihom \;\; \phiinhom]$ and $\Psi = M_h \phiinhom$. For this choice, we have $Q=I$. Note that this choice also guarantees the invertibility of $L_r$.

The proof of the theorem provides an explanation why eliminating the pressure in the velocity-only ROM is possible. 
As can be seen from
% From 
ODEs \eqref{eq:vp rom momentum eq} and \eqref{eq:vp rom hom ode}, 
% we know that 
the pressure term does not affect the coefficients of the homogeneous modes ($a_1(t)$) but only those of the inhomogeneous modes ($a_2(t)$).
% we know that the pressure term only affects the coefficients of the inhomogeneous modes of $\Phi$. 
Since these inhomogeneous coefficients are uniquely determined by the mass equation \eqref{eq:vp rom mass eq}, 
% we do not need to simulate the 
they do not require simulating
ODE \eqref{eq:vp rom momentum eq},
% for these coefficients, 
so the pressure is superfluous.
% not required.
% at all.
This observation is in line with the common understanding that the pressure is a (vector of) Lagrange multiplier(s) enforcing the velocity to satisfy the mass equation. 
This understanding also implies that the pressure does not affect the velocity beyond enforcing the mass equation. Indeed, the error between the approximated pressure $p_r(t) = \Psi b(t)$ and the FOM pressure $p_h(t)$ does not introduce an additional error in the ROM velocity.

In practice, the velocity-pressure formulation is slightly more expensive than the velocity-only formulation
because the ODE system is integrated in time for $R_V = \Rhom + \Rinhom$ instead of only $\Rhom$ coefficients and additionally the pressure coefficients have to be computed by solving the reduced Poisson equation \eqref{eq:vp rom poisson eq} at every time step.
Hence, the ROM described in Section \ref{sec:velo only rom} remains the method of choice.

Analogously to \cite{altmann2015finite}, we can decompose the velocity-pressure ROM velocity into a homogeneous and an inhomogeneous component using a $QR$ decomposition for arbitrary choices of $\Phi$ and $\Psi$ provided that $L_r$ is invertible. Provided the discretization satisfies the gradient-divergence duality \eqref{eq:spacedisc gradient divergence}, this decomposition can be formulated as a velocity-only ROM. Due to the stability problems of the QR decomposition described in \cite{altmann2015finite}, this velocity-only ROM would in general not be useful in practice.

\section{Numerical experiments}
% \label{sec: num exps mass cons}
\label{sec:num ex}

In our numerical experiments, we investigate our proposed velocity-only ROM regarding convergence behaviour, mass conservation and energy consistency. Furthermore, we verify the equivalence to the velocity-pressure ROM described in Section \ref{sec:eq to the vo rom} and analyse the runtime.

\subsection{Setup}
\subsubsection{Testcases}

\begin{figure}%[t]
	\centering
% \documentclass[tikz, margin=5mm]{standalone}
% \usepackage{tikz}
% \usetikzlibrary{shapes.geometric}
% \usetikzlibrary{positioning}
% \usetikzlibrary{positioning,shapes.multipart, fit,backgrounds,calc}

% %\tikzset{dot/.style={fill=black,circle}}
% %\tikzset{edot/.style={fill=black,circle}}

% \begin{document}
\begin{tikzpicture}[thick,black]
\tikzset{edot/.style={fill=black,diamond}}

\tikzset{triangle/.style = {fill=black!0, regular polygon, regular polygon 
sides=3}}
\tikzset{grey_triangle/.style = {fill=black!30, regular polygon, regular 
polygon sides=3}}
\tikzset{border rotated/.style = {shape border rotate=270}}
\tikzset{square/.style = {fill=black!100, rectangle}}

\draw [line width=2mm] (0,-2) -- (0,2);
\draw (10,-2) -- (10,2);
\draw (0,-2) -- (10,-2);
\draw (0,2) -- (10,2);

\draw [line width=2mm] (2,-.5) -- (2,.5);

\draw [->] (2.4,.5) -- (1.6,.5);
\draw [->] (2.4,.3) -- (1.6,.3);
\draw [->] (2.4,.1) -- (1.6,.1);
\draw [->] (2.4,-.1) -- (1.6,-.1);
\draw [->] (2.4,-.3) -- (1.6,-.3);
\draw [->] (2.4,-.5) -- (1.6,-.5);

\node at (2.6,.3) {f};
\node at (2,.9) {actuator disk};

\node [rotate = 90] at (-.5,0) {inflow boundary};

\draw [->] (0,-.7) -- (1,-.3);

\node at (.8,-.9) {$u_b(y,t)$};
\node at (.8,-1.4) {$v_b(y,t)$};

\end{tikzpicture}
% \end{document}
\caption{Sketch of testcase setup with actuator disk and time- and space-dependent boundary conditions.}
\label{fig:actuator setup}
\end{figure}
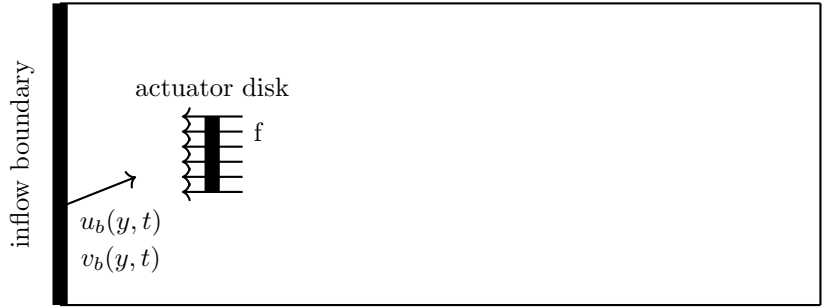

We consider two testcases based on 
% \cite[Section 6.3]{sanderse2020non}. 
\cite[Section 6.4.2]{sanderse2013energy} modelling air flow around a wind turbine. The simulation domain is $[0,10]\times [-2,2]$ with an actuator disk of length 1 at $x = 2$ as depicted in Fig.\ \ref{fig:actuator setup}. On 
the top ($y=y_\mathrm{max} = 2$), on the bottom ($y=y_\mathrm{min} = -2$) and on the right ($x=10$), we 
% define 
prescribe
%outflow 
the outflow boundary conditions
\begin{align}
	(-pI+\nu \nabla u)\cdot n = -p_\infty I \cdot n
	\label{eq:outflow bc 11}
\end{align}
with $p_\infty = 0$ and $\nu = 10^{-2}\frac{m^2}{s}$. 
% The inflow boundary on the left (x = 0) differs from \cite[Section 6.4.2]{sanderse2013energy}. 
In contrast to \cite[Section 6.4.2]{sanderse2013energy}, we consider inflow on the left boundary ($x=0$) which is not only time- but also space-dependent.
For the first testcase \consistent, we consider inflow with constant magnitude 1 but time- and space-dependent angle,
\begin{align}
    u_b(y,t) = \cos \alpha (y,t), \h 
    v_b(y,t) = \sin \alpha (y,t), \h \text{ with } \h
    \alpha(y,t) = \frac{\pi}{6}\sin(y-\frac{t}{2}).
\end{align}
As we will see in Section \ref{sec:singular values}, 
% we can approximate 
these boundary conditions 
can be approximated
well by a lower dimensional approximation $\tybc(t)$ as described in Section \ref{sec:effic rhs}. In contrast, the second testcase exemplifies the extreme case that the boundary conditions cannot be reduced effectively.
For this second testcase \movemode, we move the parabolic inflow field 
% \begin{align}
$
    u_\mathrm{parabolic}(y) = \frac{1}{10} (y-y_\mathrm{min})(y_\mathrm{max}-y)
    $
% \end{align}
uniformly such that it enters the simulation domain from the bottom at $\tstart$ and leaves it at the top at $\tend$,
\begin{align}
% $
u_b(y,t) = u_\mathrm{parabolic}\left(y + \frac{t-\tend}{\tend-\tstart}(y_\mathrm{max} - y_\mathrm{min}) \right)
% $
\hs{2} \text{ for } \tstart \leq t \leq \tend.
\end{align}
% as depicted in Figure (create!!!).
This inflow does not have a tangential component, so $v_b = 0$.

For both testcases, the Poisson matrix $L_h$ is invertible, so $\bL = L_h$ in the FOM Poisson equation \eqref{eq:fom pressure possion eq} and in the definition of the lifting function \eqref{eq:def approx lifting function}.

% Besides the inflow boundary conditions, there is no difference between the two testcases.
The actuator disk represents a momentum sink and is modeled as a constant force $f = 0.25$ acting in negative $x$-direction
and we consider the time intervals $[0,4\pi]$ and $[0,20]$, respectively.

\subsubsection{Numerical methods}
We discretize the simulation domain with a uniform $200 \times 80$ staggered grid. For time integration, we use the explicit fourth-order Runge-Kutta method RK4 \cite{butcher2008numerical}
% with time step size $\Delta t=\frac{4\pi}{800}$ for the 'varying inflow' testcase and $\Delta t = \frac{20}{800}$ for the 'moving mode' testcase, for both FOM and ROM simulations.
and discretize for both testcases the respective time interval equidistantly with time step size $\Delta t=\frac{\tend}{800}$.

The initial condition is given by $V_h^0 = \Vhinhom^0$ which is computed according to \eqref{eq:def lifting function} based on the respective boundary conditions.
For the \movemode\ testcase, this definition results in $V_h^0 = 0$ while the initial condition for the \consistent\ testcase is a nonzero function.

The POD bases $\phihom$ and $\phibc$ are constructed from the snapshot matrices
% $\Xhom = [\Vhhom^0 \;\; \dots \;\; \Vhhom^j \;\; \dots \;\; \Vhhom^{800}]$ where $\Vhhom^j$ is the time discrete approximation of $\Vhhom(t^j)$ with $t^j = j\cdot \Delta t$ and
% $X_{bc} = [\ybc(t^0) \;\;\dots\;\; \ybc(t^j) \;\;\dots\;\; \ybc(t^{800})]$.
\begin{align}
% $$
    \Xhom = [\Vhhom^0 \;\; \dots \;\; \Vhhom^j \;\; \dots \;\; \Vhhom^{800}] 
    \hs{2} \text{ and } \hs{2}
    X_{bc} = [\ybc(t^0) \;\;\dots\;\; \ybc(t^j) \;\;\dots\;\; \ybc(t^{800})],
% $$
\end{align}
where $\Vhhom^j$ is the time discrete approximation of $\Vhhom(t^j)$ with $t^j = j\cdot \Delta t$.

We compute the coefficient vector $\abc(t) = \phibc\phibc^T \ybc(t)$ 
during the offline precomputation phase for all required time steps, 
% of the boundary condition approximation $\ybc(t)$, we use the exact projection of the original boundary condition vector $\ybc(t)$ onto the POD basis $\phibc$,
as described in Section \ref{sec:effic rhs}.

The numbers of ROM modes $\Rhom$ and $\Rbc$ can be chosen independently of each other. In the numerical experiments here, however, we choose $\Rhom = \Rbc = R$ for convenience only, and consider velocity-only ROMs for different values of $R$.
% We consider velocity-only ROMs with $\Rhom = \Rbc = R$ modes.

% For both testcases, we investigate velocity-only ROMs for all three discussed choices of the bases $\phibc$, $\phiinhom$, $\phihom$. 
% The first two bases choices \mthesis and \optimal consist of $\phibc$ as described in Section \ref{sec:approximation_bc} and Section \label{sec:velo consis bc approx}, respectively, $\phiinhom$ as described in Section \ref{sec:phi_inhom} and $\phihom$ as described in Section \ref{sec:hhd inspired decomposition}. The third bases choice \qr consists of $\phibc$ as described in Section \ref{sec:vp rom} and $\phiinhom$ and $\phihom$ as described in Section \ref{sec:equivalence to vo rom}.
% For all three velocity-only ROMs resulting from these bases choices we construct corresponding velocity-pressure ROMs with $\phi = [\phihom \;\; \phiinhom]$, $\psi$ obtained from orthonormalization of $M_h\phiinhom$ and $\phibc$ as in the corresponding velocity-only ROM.

% To avoid machine precision errors,
% we initialize the velocity-pressure ROMs with the same coefficient vectors as the velocity-only ROMs.

\subsubsection{Singular value decays}
\label{sec:singular values}
\begin{figure}%[t]
	\begin{subfigure}{0.5\textwidth}
	\centering
\includegraphics[width=\textwidth]{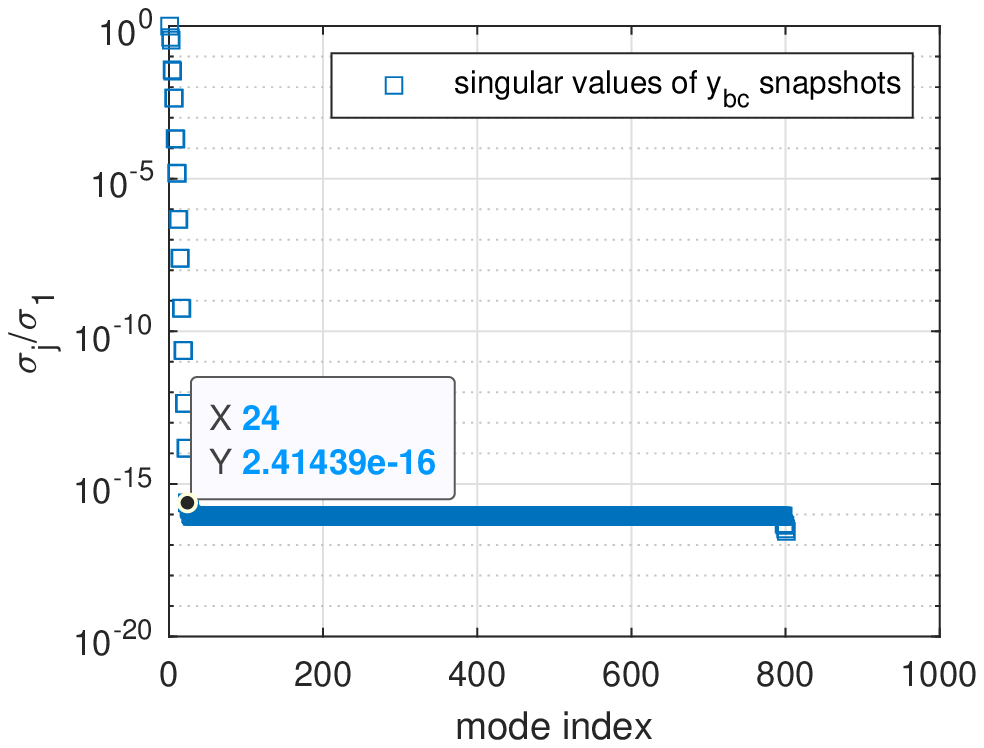}
\caption{\Consistent}
\label{fig:singular vals ybc consistent}
	\end{subfigure}
	\begin{subfigure}{0.5\textwidth}
	\centering
\includegraphics[width=\textwidth]
{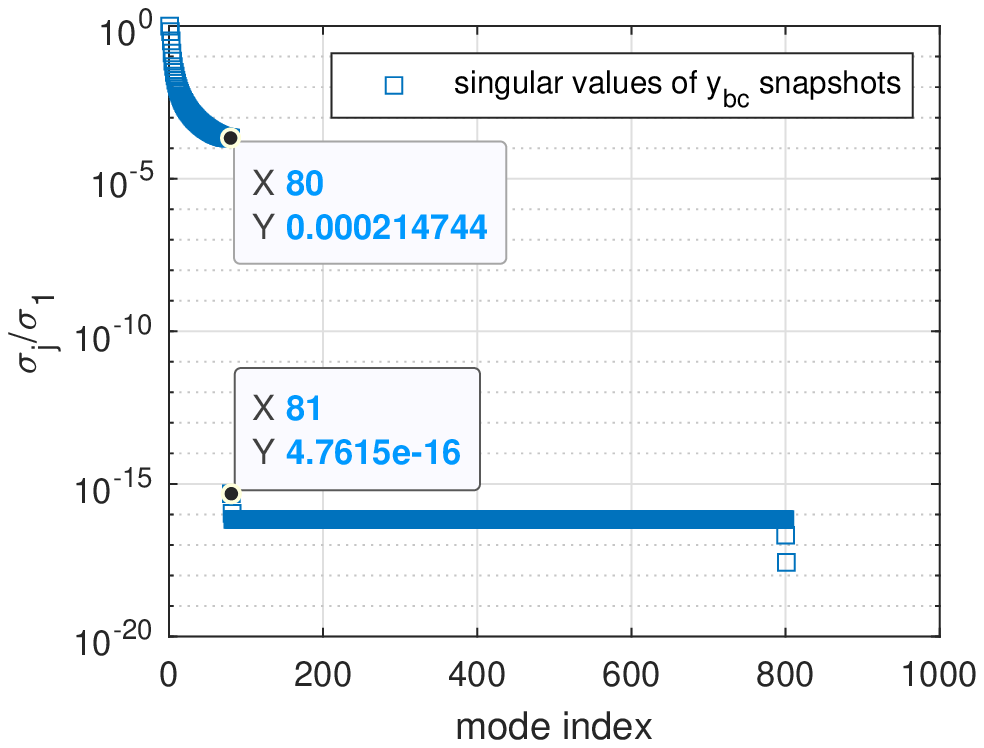}
\caption{\Movemode}
\label{fig:singular vals ybc movemode}
	\end{subfigure}
\caption{Singular values of the snapshot matrices
% $X_{bc} = [\ybc(t^0) \;\;\dots\;\; \ybc(t^j) \;\;\dots\;\; \ybc(t^{800})]$ 
$\Xbc$
for the two testcases divided by largest singular value.}
\label{fig:singular vals ybc}
\end{figure}

First, we investigate the singular value decays of the snapshot matrices which we construct our ROM bases from.
As can be seen in Fig.\ \ref{fig:singular vals ybc consistent}, only a few singular values of the boundary condition vector snapshot matrix $\Xbc$ have a significant magnitude for the \consistent\ testcase. 
% In fact,
% only 23 of them 
% % are at least $10^{-15}$ times as large as the largest singular value. 
% are larger than machine precision, when divided by the largest singular value.

In contrast, Fig.\ \ref{fig:singular vals ybc movemode} shows that 80 singular values of $\Xbc$ for the \movemode\ testcase have a relative magnitude larger than $10^{-4}$. This number equals the number of 
% nonzero entries of $y_M$, 
grid points at the inflow boundary and the number of nonzero entries in $y_M$,
so these entries do not exhibit any lower-dimensional coherent structure.

% \begin{figure}%[t]
% 	\begin{subfigure}{0.5\textwidth}
% 	\centering
% \includegraphics[width=\textwidth]
% {pics2/singular_values_yM_snapshots_consistent.eps}
% \caption{Varying inflow.}
% \label{fig:singular vals yM consistent}
% 	\end{subfigure}
% 	\begin{subfigure}{0.5\textwidth}
% 	\centering
% \includegraphics[width=\textwidth]
% {pics2/singular_values_yM_snapshots_movemode.eps}
% \caption{Moving mode.}
% \label{fig:singular vals yM movemode}
% 	\end{subfigure}
% \caption{Singular values of the snapshot matrices $X_{bc}$ wrong !!! for the two testcases.}
% \label{fig:singular vals yM}
% \end{figure}

\begin{figure}%[t]
	\begin{subfigure}{0.5\textwidth}
	\centering
\includegraphics[width=\textwidth]
{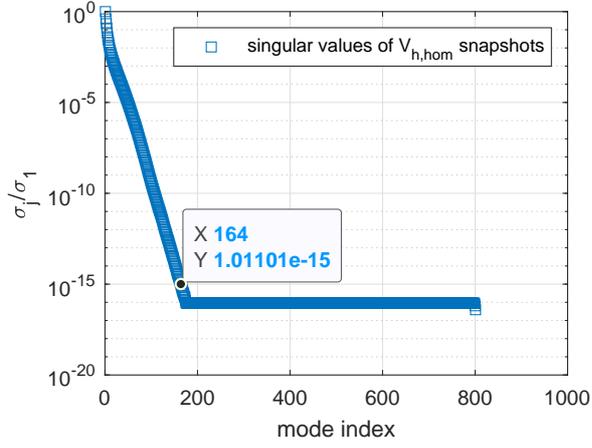}
\caption{\Consistent}
\label{fig:singular vals Vhom consistent}
	\end{subfigure}
	\begin{subfigure}{0.5\textwidth}
	\centering
\includegraphics[width=\textwidth]
{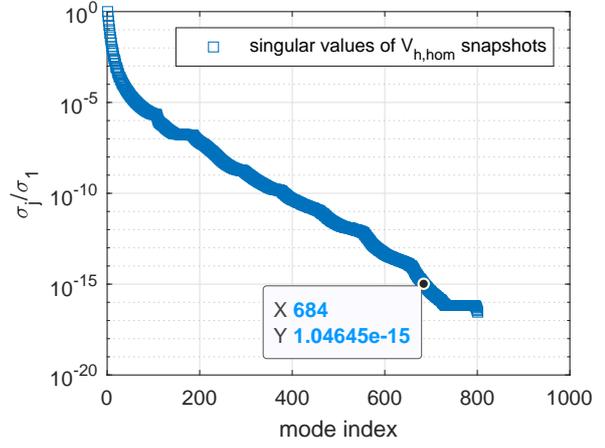}
\caption{\Movemode}
\label{fig:singular vals Vhom movemode}
	\end{subfigure}
\caption{Singular values of the snapshot matrices
% $\Xhom = [\Vhom^0 \;\; \dots \;\; \Vhom^j \;\; \dots \;\; \Vhom^{800}]$ 
$\Xhom$
for the two testcases.}
\label{fig:singular vals Vhom}
\end{figure}

Similarly, the decays of singular values of the velocity snapshot matrices $\Xhom$ shown in Fig.\ \ref{fig:singular vals Vhom} differ considerably. While only 164 singular values 
are larger than machine precision,
when divided by the largest singular value,
% are at least $10^{-15}$ times as large as the largest one 
for the \consistent\ testcase, 684 singular values are larger than this relative threshold for the \movemode\ testcase.

% \begin{figure}%[t]
% 	\begin{subfigure}{0.5\textwidth}
% 	\centering
% \includegraphics[width=\textwidth]
% {pics2/singular_values_Vh_consistent3.eps}
% \caption{Varying inflow.}
% \label{fig:singular vals Vh consistent}
% 	\end{subfigure}
% 	\begin{subfigure}{0.5\textwidth}
% 	\centering
% \includegraphics[width=\textwidth]
% {pics2/singular_values_Vh_movemode3.eps}
% \caption{Moving mode.}
% \label{fig:singular vals Vh movemode}
% 	\end{subfigure}
% \caption{Singular values of the snapshot matrices $X_{h} = [V_h^0 \;\; \dots \;\; V_h^j \;\; \dots \;\; V_h^{800}]$ for the two testcases.}
% \label{fig:singular vals Vh}
% \end{figure}

\subsection{Results}
\subsubsection{Convergence}

\begin{figure}%[t]
	\begin{subfigure}{0.5\textwidth}
	\centering
\includegraphics[width=\textwidth]
{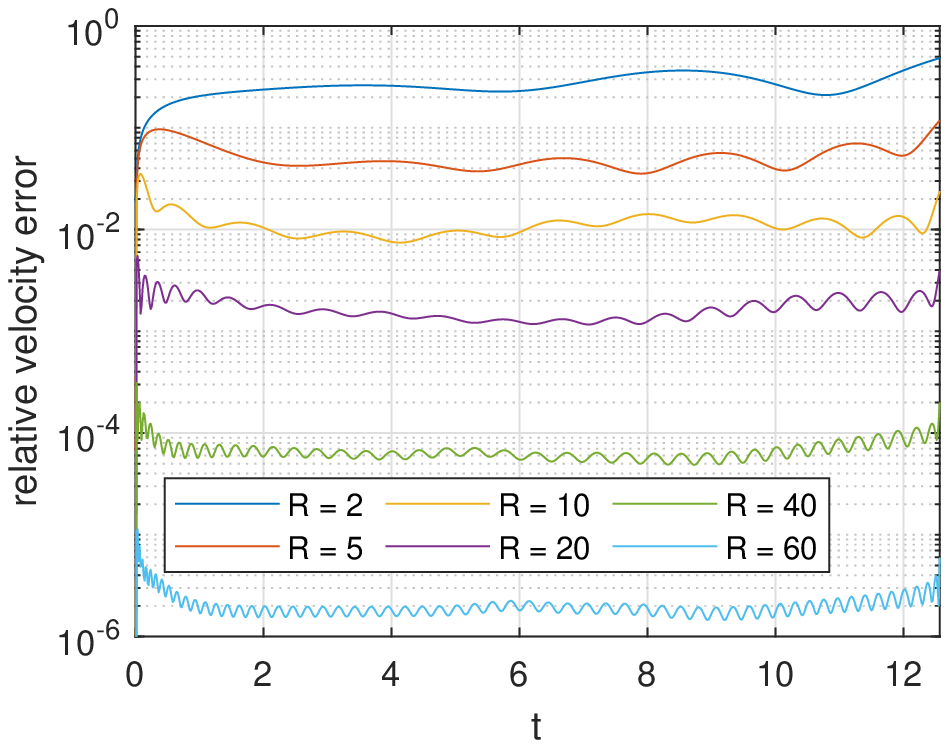}
\caption{\Consistent}
\label{fig:velocity error consistent}
	\end{subfigure}
	\begin{subfigure}{0.5\textwidth}
	\centering
\includegraphics[width=\textwidth]
{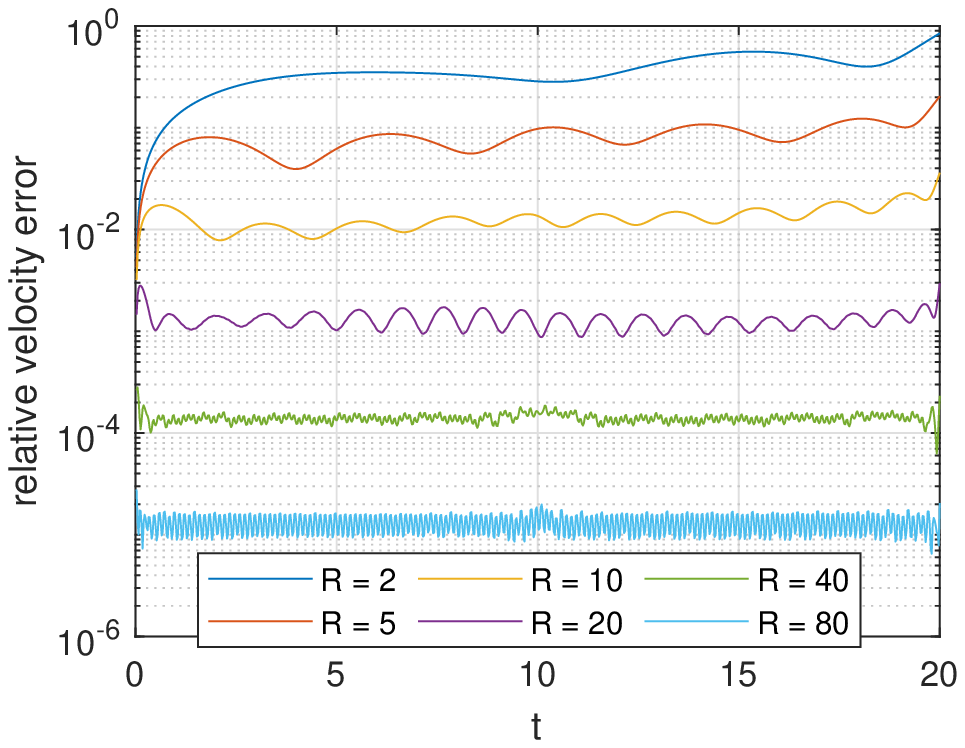}
\caption{\Movemode}
\label{fig:velocity error movemode}
	\end{subfigure}
\caption{
Relative velocity error
$\frac{\|V_h^j - V_r^j\|_\Oh}{\average_j\left[\|V_h^j\|_\Oh\right]}$
for the two testcases. 
% $\Oh$-norm of difference between FOM and ROM velocity divided by $\Oh$-norm of FOM velocity.
% ROM velocity error and velocity best approximation error. $\Oh$-norm of difference to the FOM velocity
% % Singular values of the snapshot matrices $X_{h}$ 
% for the two testcases.
}
\label{fig:velocity error}
\end{figure}

First, we investigate whether the ROM converges to the FOM for increasing number of modes. 
% For this purpose, we consider the relative velocity error 
% % $\frac{\|V_h^j - V_r^j\|_\Oh}{\|V_h^j\|_\Oh}$ 
% of ROMs with $R = \Rhom = \Rbc$ modes.
As can be seen in Fig.\ \ref{fig:velocity error} the relative velocity error clearly decreases monotonously for increasing number of ROM modes.
Consequently, our velocity-only ROM for time-dependent boundary conditions indeed provides a convergent approximation to the FOM.

% For the same simulations as in the previous subsection, Fig.\ \ref{fig:velocity error} shows the ROM velocity errors. For both testcases, we clearly observe convergence upon increasing the number of ROM modes. 

\subsubsection{Mass conservation}

\begin{figure}%[t]
	\begin{subfigure}{0.5\textwidth}
	\centering
\includegraphics[width=\textwidth]
{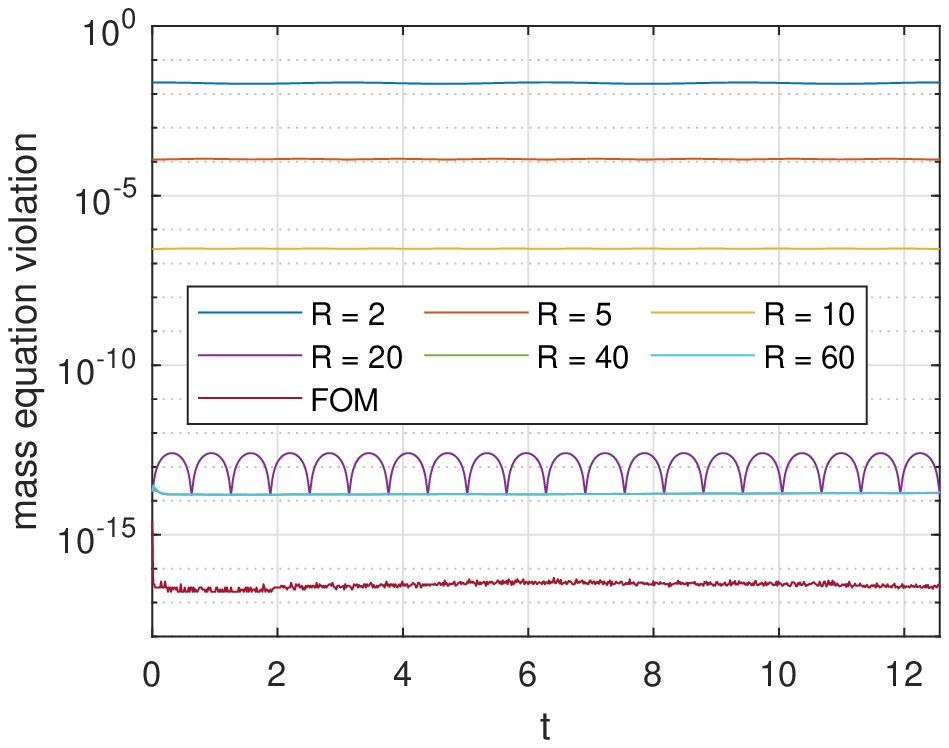}
\caption{\Consistent}
\label{fig:mass viola consistent}
	\end{subfigure}
	\begin{subfigure}{0.5\textwidth}
	\centering
\includegraphics[width=\textwidth]
{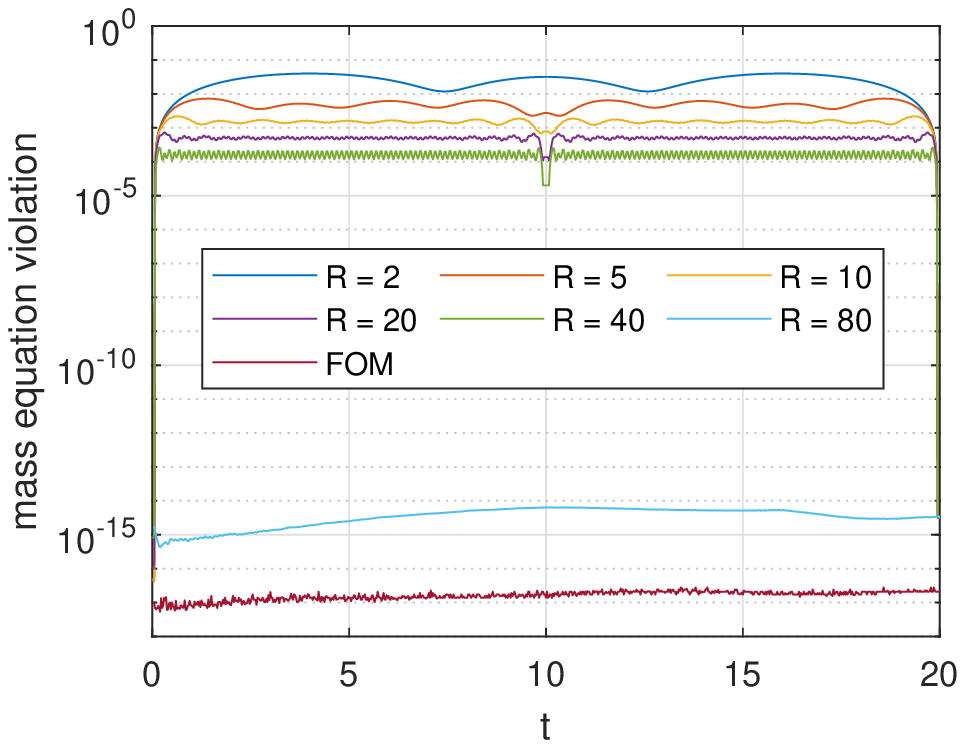}
\caption{\Movemode}
\label{fig:mass viola movemode}
	\end{subfigure}
\caption{Mass conservation violation
$\|M_h V_r^j - y_M(t^j)\|_2$ for the two testcases.
% $L^2$-norm of the residual of the mass equation with respect to the original boundary conditions
% % Singular values of the snapshot matrices $X_{h}$ 
% for the two testcases.
}
\label{fig:mass viola}
\end{figure}

Second, we investigate to which extent our proposed ROM violates 
% the mass equation with respect to the original boundary conditions 
the original mass equation
\eqref{eq:rom exact mass eq}. 
% This violation is shown in Fig.\ \ref{fig:mass viola} for the same simulations as considered in the previous section.
As explained in Section \ref{sec:effic rhs}, the proposed ROM satisfies the mass equation exactly,  but with respect to the approximated boundary conditions. The difference between the two is shown in Fig.\ \ref{fig:mass viola}, effectively indicating how well the mass equation right-hand side $y_M(t)$ is approximated by $\tilde y_M(t)$. 
% = F_M\phibc \abc(t)$.

Indeed, as can be seen in Fig.\ \ref{fig:mass viola} the mass conservation violation decreases for increasing number of boundary condition approximation modes. 
For the testcase \consistent, the mass conservation violation is in the order of machine precision for $R\geq 20$. 
In contrast, for the \movemode\ testcase, 
% we require 
$R=80$
modes are required
to achieve a similar mass conservation accuracy. Even for $R=40$ modes, the mass conservation violation is larger than $10^{-4}$. This observation can directly be explained by the singular value decay depicted in Fig. \ref{fig:singular vals ybc movemode},  which shows that the 80 largest singular values all have a relative magnitude larger than $10^{-4}$. Note however, that in this case the boundary conditions were chosen on purpose to be highly `irreducible' with POD and as such this can be seen as a limiting case. Furthermore, the violation of the exact mass equation seems not to impair the accuracy of the actual solution $V_r$, as can be observed by comparing Fig.\ \ref{fig:velocity error consistent} to \ref{fig:velocity error movemode}.

% In Fig.\ \ref{fig:mass viola consistent}, we can also see the implications of including POD modes whose singular values have small relative magnitude. Indeed, we identify an oscillating curve for $R=20$. These

% For the testcase \consistent, $\Mbc = 40$ is for both bases choices sufficient to achieve mass conservation up to machine precision. For the testcase \movemode, the \mthesis ROM achieves mass conservation up to machine precision for $\Mbc = 80$ whereas the \optimal ROM violates the mass conservation even for $\Mbc = 80$ by more than $10^{-6}$ for many time steps.

\subsubsection{Energy consistency}

Third, we investigate whether the structural correspondence of the FOM and ROM kinetic energy evolutions derived in Section \ref{sec:kin energy evol} is accompanied by small energy errors in practice.

\begin{figure}%[t]
	\begin{subfigure}{0.5\textwidth}
	\centering
\includegraphics[width=\textwidth]
{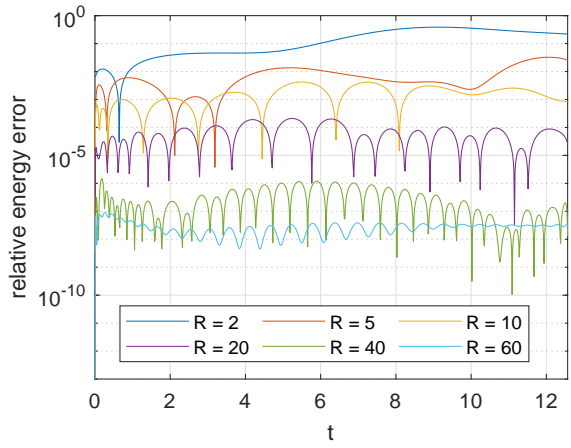}
\caption{\Consistent}
\label{fig:kin en error consistent}
	\end{subfigure}
	\begin{subfigure}{0.5\textwidth}
	\centering
\includegraphics[width=\textwidth]
{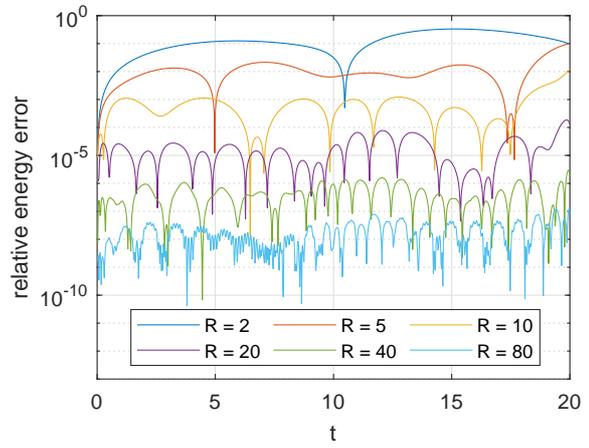}
\caption{\Movemode}
\label{fig:kin en error movemode}
	\end{subfigure}
\caption{Relative kinetic energy error $\frac{\left|K_h^j-K_r^j\right|}{\average_j \left[K_h^j\right]}$ for both testcases.
}
\label{fig:kin en error}
\end{figure}

\begin{figure}%[t]
	\begin{subfigure}{0.5\textwidth}
	\centering
\includegraphics[width=\textwidth]
{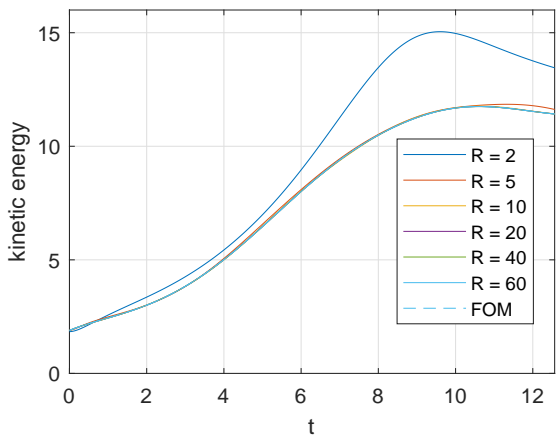}
\caption{\Consistent}
\label{fig:kin energy consistent}
	\end{subfigure}
	\begin{subfigure}{0.5\textwidth}
	\centering
\includegraphics[width=\textwidth]
{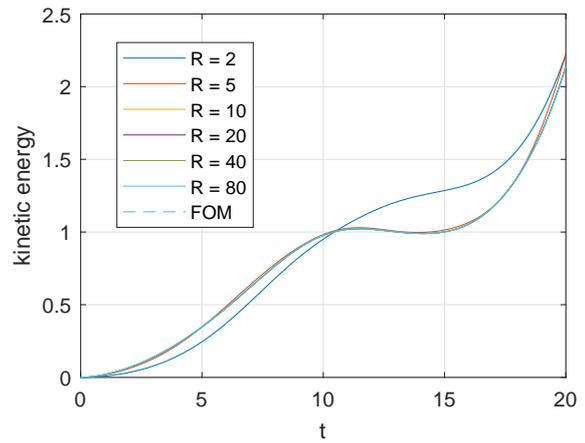}
\caption{\Movemode}
\label{fig:kin energy movemode}
	\end{subfigure}
\caption{Time evolution of the kinetic energy for both testcases.
}
\label{fig:kin energy}
\end{figure}

Indeed, Fig.\ \ref{fig:kin en error} shows that the relative kinetic energy errors are approximately at the same order of magnitude as the corresponding relative velocity errors depicted in Fig.\ \ref{fig:velocity error}, and the behavior of the total FOM kinetic energy in time is accurately reproduced, see Fig.\ \ref{fig:kin energy}. 
% We stress that the kinetic energy is \textit{not} constant in time for the considered testcases. our ROMs do not merely keep the energy constant but accurately reproduce the energy variation of the FOM.

\subsubsection{Verification of ROM equivalence}

\begin{figure}%[t]
	\begin{subfigure}{0.5\textwidth}
	\centering
\includegraphics[width=\textwidth]
{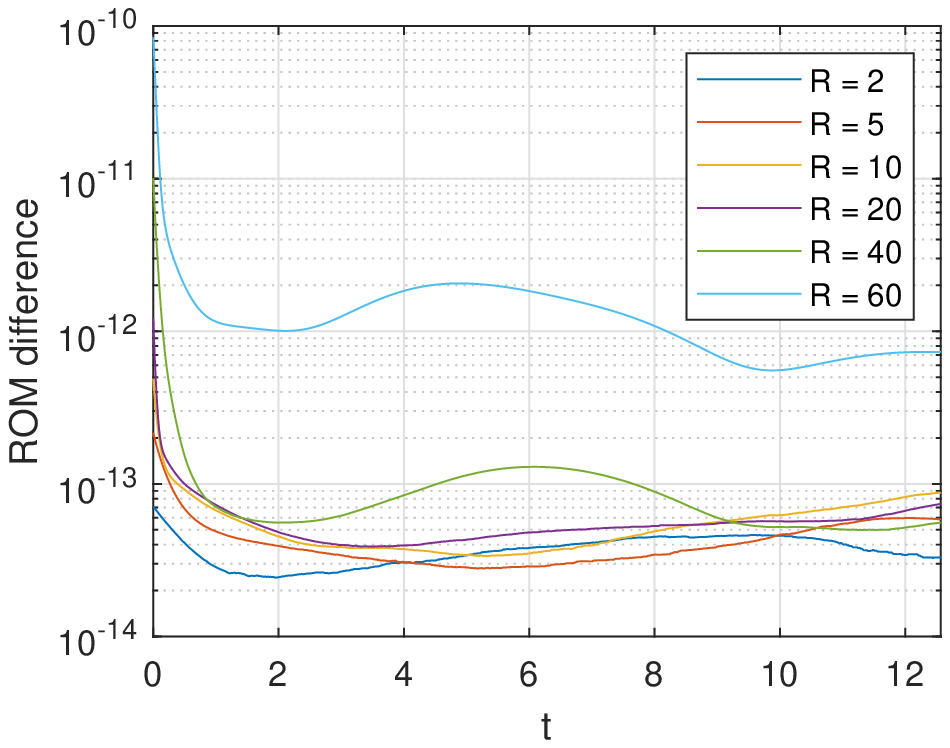}
\caption{\Consistent}
\label{fig:rom equivalence consistent}
	\end{subfigure}
	\begin{subfigure}{0.5\textwidth}
	\centering
\includegraphics[width=\textwidth]
{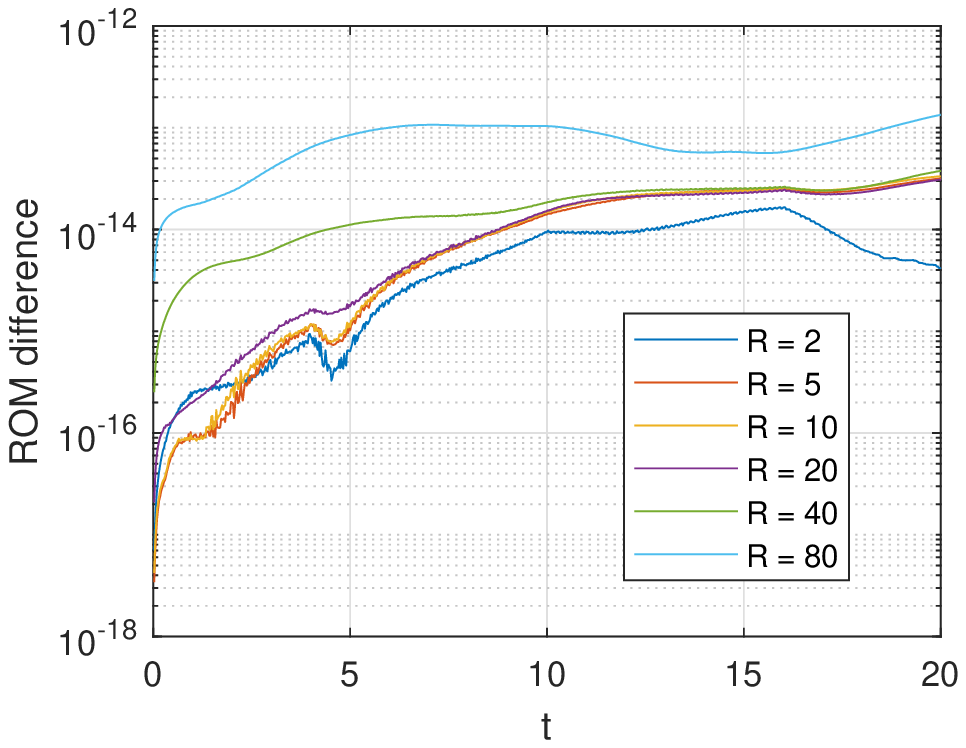}
\caption{\Movemode}
\label{fig:rom equivalence movemode}
	\end{subfigure}
\caption{ROM equivalence error
$\left\|\Vrvo^j - \Vrvp^j\right\|_\Oh$ for both testcases.
% . $\Oh$-norm of difference of velocities of velocity-only ROMs with $\Rhom = \Rinhom = \Rbc = R$ and velocity-pressure ROMs with $R_V = 2R$, $R_p = \Rbc = R$.
}
\label{fig:rom equivalence}
\end{figure}

Next, we verify the 
% To verify the 
ROM equivalence demonstrated in Section 
% \ref{sec:equivalence to vo rom}.
\ref{sec:eq to the vo rom}.
% , we compute the difference between the velocities of the velocity-only and the velocity-pressure ROM. 
For this purpose, we simulate the same testcases as in the previous experiment but now with the velocity-pressure ROM \eqref{eq:vp rom mass eq}-\eqref{eq:vp rom momentum eq} with the bases $\Phi = [\phihom \;\;\phiinhom]$ and $\Psi = M_h\phiinhom$ as described in Section \ref{sec:eq to the vo rom}. Here, we compute $\phiinhom$ via a QR decomposition of $\Finhom$. Hence, the numbers of modes of the velocity-pressure ROM are given by $\Rbc = \Rinhom = R_p = R$ and $R_V = \Rhom + \Rinhom = 2R$.
As can be seen in Fig.\ \ref{fig:rom equivalence}, the error between the two ROMs ranges for all values of $R$ between $10^{-18}$ and $10^{-10}$. Hence, the ROM equivalence is confirmed up to machine precision.

Note that we observe larger errors for ROMs with larger number of modes, which is likely because these ROMs accumulate larger numbers of machine precision errors.
% Furthermore, we observe the errors of the \consistent testcase in Fig.\ \ref{fig:rom equivalence consistent} to be consistently

\subsubsection{Runtime analysis}
\begin{figure}%[t]
	\centering
\includegraphics[width=.6\textwidth]
{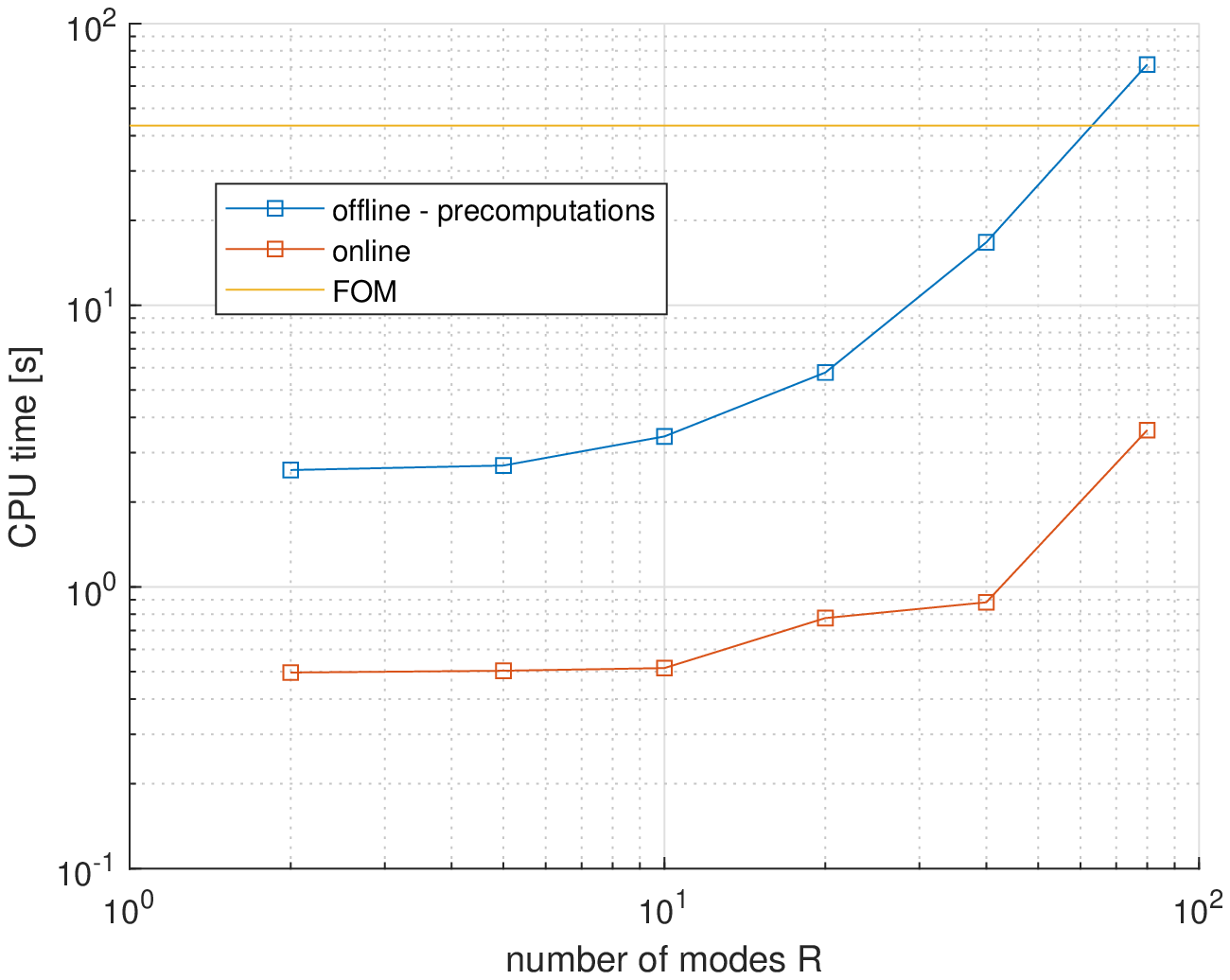}
\caption{Computation time for testcase \movemode\ with $\Rhom = \Rinhom = \Rbc = R$.}
\label{fig:computation time}
\end{figure}

Fig.\ \ref{fig:computation time} shows the computation time of the ROM in comparison with the FOM. For all considered numbers of modes, we observe a significant speed-up of the ROM online phase compared to the FOM simulations. As in \cite{sanderse2020non}, this speed-up amounts to approximately two orders of magnitude for moderate numbers of modes.
The offline precomputations, on the other hand, have significant costs and even exceed the FOM simulation costs for $R = 80$.

However, when used in practice for a parametric study, the ROM would already pay off for a small number of ROM simulations.

% Since these precomputations are performed only once, the use of the ROM pays off also for a small number of ROM simulations.

\section{Conclusion}
\label{sec: results mass cons}
In this paper we have generalized the 
% mass-conserving 
energy-conserving
ROM for the incompressible Navier-Stokes equations proposed by Sanderse \cite{sanderse2020non} to arbitrary time-dependent boundary conditions. 
% Like 
% Sanderse's ROM, our proposed ROM is velocity-only. Hence, the computation of the 
% pressure is not required and as a consequence, pressure-related $\inf$-$\sup$ instability is avoided. 
% As in \cite{sanderse2020non}, 
Our ROM is velocity-only, which avoids computing the pressure and pressure-related $\inf$-$\sup$ instability, and it preserves the structure of the kinetic energy evolution.
The ROM is based on a Helmholtz-Hodge-inspired decomposition that separates the velocity into two orthogonal components: i) a velocity that satisfies the mass equation with homogeneous boundary conditions and ii) a lifting function that satisfies the mass equation with the prescribed time-dependent boundary conditions. 
% At any point in time, the lifting function is defined by the mass equation right-hand side at this point in time only. 
The homogeneous component is defined by the Galerkin-projected momentum equation. Due to the divergence-free projection basis, 
% we can eliminate the pressure term in the momentum equation to find a velocity-only ODE.
the pressure term disappears form the momentum equation resulting in the proposed velocity-only ROM.

To allow efficient simulation of the ROM, we approximate the prescribed boundary conditions by a lower-dimensional approximation via POD. The resulting ROM satisfies the mass equation with respect to these approximated boundary conditions. 
Our numerical experiments show that the error introduced by this boundary condition approximation does not dominate the velocity error. In these experiments, we use as many modes to approximate the boundary conditions as we use to approximate the homogeneous velocity field. In future research, different choices of these mode numbers should be investigated as well.

% We have derived a kinetic energy evolution on the ROM level that matches the one on the FOM level. 
We have theoretically shown that the orthogonality of the components leads to a ROM kinetic energy evolution that has the same structure as the FOM kinetic energy evolution.
% In future work, this finding can be used to infer nonlinear stability of the proposed ROM for specific boundary conditions.
Our numerical experiments confirm that the proposed ROM is able to accurately mimic the kinetic energy evolution.
In future work, the structural equivalence of the kinetic energy evolution can be used to derive bounds on the kinetic energy and hence to derive nonlinear stability results.

Furthermore,
we have proven the important result that the proposed velocity-only ROM is equivalent to a velocity-pressure ROM, when employing a specific choice of basis functions, and 
% used this equivalence to derive a kinetic energy evolution on the ROM level that matches the one on the FOM level. 
confirmed this equivalence in our numerical experiments.
This equivalence can function as a key concept to understand the role of the pressure in ROMs for the incompressible Navier-Stokes equations. It gives important insights into the behaviour of the velocity-pressure ROM regarding the choice for the number of velocity, pressure and boundary condition modes.
% In addition, the equivalence of existing velocity-pressure ROMs to velocity-only ROMs can help understanding their behaviour.
% regarding changes in the numbers of velocity and pressure modes.
While this ROM equivalence is based on the assumption of the discrete gradient-divergence duality, future research should investigate whether discretizations that do not exhibit this property also allow for a similar equivalence statement.

In numerical experiments, we have shown the convergence of the ROM upon increasing the number of ROM modes and the accuracy of the boundary condition approximation. 
% Furthermore, we have confirmed the equivalence between the velocity-only and the velocity-pressure ROM.
% In our numerical experiments, we have been able to show that the proposed mass-conserving ROMs violate the mass conservation equation indeed only in the range of machine precision. Furthermore, in the second test case, we observed the ROM errors to be significantly smaller for mass-conserving ROMs than the for standard POD-Galerkin ROMs with the same number of POD modes. 
Instead of performing a hyper-reduction technique
for the convective term,
we have used an exact offline decomposition.
% , we proposed an offline decomposition that introduces a consistent approximation of the boundary conditions and provides the exact computation of the ROM ODE RHS with respect to this boundary condition approximation. 
Scaling cubically in the number of modes, this offline decomposition is limited to small numbers of modes. For moderate numbers of modes, the speed-up of the ROM online phase compared to the FOM is approximately two orders of magnitude.
If larger numbers of modes are required that prohibit this cubic scaling, we suggest to incorporate energy-conserving hyper-reduction methods \cite{klein2022structure}.

\section*{CRediT authorship contribution statement}
\textbf{H. Rosenberger:} Conceptualization, Methodology, Software, Writing - original draft.
\textbf{B. Sanderse:} Conceptualization, Methodology, Software, Writing - review \& editing, Funding acquisition.

\section*{Acknowledgements}
This publication is part of the project "Discretize first, reduce next" (with project number VI.Vidi.193.105 of the research programme NWO Talent Programme Vidi which is (partly) financed by the Dutch Research Council (NWO).

\FloatBarrier

% Appendix
% \appendix
% \input{appendix/discretization_of_bc}
% \input{sections/skew_symmetry_summary}
% \input{appendix/velocity-only_pressure_recovery}
% \input{appendix/velo-consistent_BC_approx}
% \input{appendix/boundary_energy}
% \input{appendix/pres_ROM_ini}
% \input{appendix/kin_ener_vo_rom}

% \input{sections/mindblow}
% \input{sections/consistent_bases}
% \input{sections/equivalence_theorems}

% \input{appendix/standard_ROM}

% \input{sections/kinetic_energy_evolution}
% \input{sections/skew-symmetry}
% \input{sections/dong_bc}
% \input{sections/botchOBC}
% \input{sections/dirichlet_bound}
% \input{sections/skew_symmetry_summary}
% \input{sections/foias}
% \input{sections/lions}
% \input{sections/lions_conti}
% \input{sections/POD_ROMs}

% \input{section_archive/kinetic_energy_evolution}
% \input{section_archive/skew-symmetry}
% \input{section_archive/dirichlet_bound}
% \input{section_archive/skew_symmetry_summary}
% \input{section_archive/POD_ROMs}

% \input{appendix/boundary_energy}
% \input{section_archive/foias}

% \input{sections/formulas}

% References
% \section*{References}
\bibliographystyle{elsarticle-num} 
\bibliography{reference}

\end{document}